\newif\ifpdf
\ifx\pdfoutput\undefined
\pdffalse 
\else
\pdfoutput=1 
\pdftrue \fi

\newif\ifarxiv
\arxivtrue

\documentclass[probth,final]{svjour}
\usepackage{cite}
\usepackage{amsmath}
\usepackage{amsfonts}
\usepackage{amssymb}
\usepackage{epsf}
\usepackage{graphicx}
\usepackage{graphics}
\usepackage{verbatim}
\usepackage{epsfig}

\IfFileExists{myowntimes.sty}
	{\usepackage{myowntimes}}
	{\usepackage{times}\usepackage{mathrsfs}}
	
\DeclareFontFamily{OT1}{fraktura}{}
\DeclareFontShape{OT1}{fraktura}{m}{n} {<5> <6> <7> <8> <9> <10> <11> <12> <13> <14.4> [1.1] eufm10}{}
\DeclareMathAlphabet{\fraktura}{OT1}{fraktura}{m}{n}

\ifarxiv
\fi

\newenvironment{proofsect}[1]{\vskip0.1cm\noindent{\rmfamily\itshape #1.}}{\qed\vspace{0.15cm}}

\spnewtheorem{mylemma}[theorem]{Lemma}{\bf}{\it}
\spnewtheorem{myproposition}[theorem]{Proposition}{\bf}{\it}
\spnewtheorem{mycorollary}[theorem]{Corollary}{\bf}{\it}
\spnewtheorem{mydefinition}[theorem]{Definition}{\bf}{\it}
\spnewtheorem{myquestion}{Question}{\bf}{\it}
\spnewtheorem{myconjecture}[myquestion]{Conjecture}{\bf}{\it}
\spnewtheorem{myremark}[theorem]{Remark}{\it}{\rm}
\numberwithin{equation}{section}
\numberwithin{theorem}{section}

\newcommand{\supp}{\operatorname{supp}}

\newcommand{\textd}{\text{\rm d}\mkern0.5mu}
\newcommand{\texti}{\text{\rm  i}\mkern0.7mu}
\newcommand{\texte}{\text{\rm  e}\mkern0.7mu}

\newcommand{\mini}{{\text{\rm min}}}

\renewcommand{\mod}{\,\text{\rm mod}\,}

\renewcommand{\AA}{\mathcal A}
\newcommand{\BB}{\mathcal B}

\newcommand{\DD}{\mathcal D}

\newcommand{\GG}{\mathcal G}

\newcommand{\TT}{\mathcal T}

\newcommand{\XX}{\mathcal X}

\newcommand{\B}{\mathbb B}

\newcommand{\E}{\mathbb E}

\newcommand{\BbbP}{\mathbb P}

\newcommand{\R}{\mathbb R}

\newcommand{\T}{\mathbb T}

\newcommand{\Z}{\mathbb Z}

\newcommand{\scrF}{\mathscr{F}}
\newcommand{\scrG}{\mathscr{G}}

\newcommand{\scrT}{\mathscr{T}}

\newcommand{\twoeqref}[2]{(\ref{#1}--\ref{#2})}
\newcommand{\1}{{1\mkern-4.5mu\textrm{l}}}
\renewcommand{\1}{\text{\sf 1}}

\newcommand{\bk}{\boldsymbol k}
\newcommand{\bp}{\boldsymbol p}
\newcommand{\bq}{\boldsymbol q}
\newcommand{\bzero}{\boldsymbol 0}

\def\myffrac#1#2 in #3{\raise 2.6pt\hbox{$#3 #1$}\mkern-1.5mu\raise 0.8pt\hbox{$#3/$}\mkern-1.1mu\lower 1.5pt\hbox{$#3 #2$}}
\newcommand{\ffrac}[2]{\mathchoice%
	{\myffrac{#1}{#2} in \scriptstyle}
	{\myffrac{#1}{#2} in \scriptstyle}
	{\myffrac{#1}{#2} in \scriptscriptstyle}
	{\myffrac{#1}{#2} in \scriptscriptstyle}
}

\newcommand{\Var}{\text{\rm Var}}

\newcommand{\cc}{{\text{\rm c}}}

\newcommand{\zz}{\fraktura z}

\newcommand{\wt}{\widetilde}

\newcommand{\kappaX}[1]{\mathchoice%
	{\kappa_{\text{\rm\fontsize{6.5}{6}\selectfont #1}}}
	{\kappa_{\text{\rm\fontsize{6.5}{6}\selectfont #1}}}
	{\kappa_{\text{\rm\fontsize{5}{5}\selectfont #1}}}
	{\kappa_{\text{\rm\fontsize{5}{5}\selectfont #1}}}
}
\newcommand{\kappaO}{\kappaX{O}}
\newcommand{\kappaD}{\kappaX{D}}

\newcommand{\Oind}{\text{\rm O}}
\newcommand{\Dind}{\text{\rm D}}

\newcommand{\MP}{\text{\rm MP}}
\newcommand{\MA}{\text{\rm MA}}
\newcommand{\UO}{\text{\rm UO}}
\newcommand{\UD}{\text{\rm UD}}

\def\ZO_#1{Z_{#1,\Oind}}
\def\ZD_#1{Z_{#1,\Dind}}
\def\ZMP_#1{Z_{#1,\MP}}
\def\ZMA_#1{Z_{#1,\MA}}
\def\ZUO_#1{Z_{#1,\UO}}
\def\ZUD_#1{Z_{#1,\UD}}

\newcommand{\Bvert}{\B^{\text{\rm vert}}}
\newcommand{\Bhor}{\B^{\text{\rm hor}}}
\newcommand{\Beven}{\B^{\text{\rm even}}}
\newcommand{\Bodd}{\B^{\text{\rm odd}}}

\newcommand{\ord}{\text{\rm ord}}
\newcommand{\dis}{\text{\rm dis}}
\newcommand{\pt}{p_{\text{\rm t}}}

\newcommand{\llambda}{r}
\newcommand{\hate}{\hat{\text{\rm e}}}
\newcommand{\etavert}{\eta_{\text{\rm vert}}}
\newcommand{\etahor}{\eta_{\text{\rm hor}}}

\newcommand{\wh}{\widehat}
\newcommand{\NO}{N_{\text{\rm O}}}
\newcommand{\ND}{N_{\text{\rm D}}}

\begin{document}
\title{Phase coexistence of gradient Gibbs states}
\titlerunning{Gradient Gibbs states}
\author {Marek~Biskup\inst{1} \and Roman Koteck\'y\inst{2}}
\authorrunning{M.~Biskup and R.~Koteck\'y}

\institute{Department of Mathematics, UCLA, Los Angeles, California, USA
\and
Center for Theoretical Study, Charles University, Prague, Czech Republic}

\date{}
\maketitle

\renewcommand{\thefootnote}{}
\footnotetext{\vglue-0.41cm\footnotesize\copyright\,2006 by M.~Biskup and R.~Koteck\'y. Reproduction, by any means, of the entire article for non-commercial purposes is permitted without charge.}
\renewcommand{\thefootnote}{\arabic{footnote}}

\begin{abstract}
We consider the 
(scalar) gradient fields $\eta=(\eta_b)$---with~$b$ denoting the nearest-neighbor edges in~$\Z^2$---that are distributed according to the Gibbs measure proportional to $\texte^{-\beta
H(\eta)}\nu(\textd\eta)$. Here $H=\sum_bV(\eta_b)$ is the Hamiltonian,~$V$ is a symmetric potential, ~$\beta>0$ is the inverse temperature, and~$\nu$ is the Lebesgue measure on the linear
space defined by imposing the loop condition $\eta_{b_1}+\eta_{b_2}=\eta_{b_3}+\eta_{b_4}$ for each plaquette~$(b_1,b_2,b_3,b_4)$ in~$\Z^2$.
For convex~$V$, Funaki and Spohn have shown that ergodic infinite-volume Gibbs measures are characterized by their tilt. 
We describe a mechanism by which the gradient Gibbs measures with non-convex~$V$ undergo a structural, order-disorder phase transition at some intermediate value of inverse temperature~$\beta$. 
At the transition point, there are at least two distinct gradient measures with zero tilt, 
i.e.,
$E \eta_b=0$.
\end{abstract}

\section{Introduction}
\subsection{Gradient fields}
\label{S:GGm}\noindent
One of the mathematical challenges encountered in the study of systems exhibiting phase coexistence is an
efficient description of microscopic phase boundaries. Here various levels of detail are in general possible:  
The finest level is typically associated with a statistical-mechanical model (e.g., a lattice gas) in which both the interface 
and the surrounding phases are represented microscopically; 
at the coarsest level the interface is viewed as a macroscopic (geometrical) surface between two structureless bulk phases. 
An intermediate approach  is  based on effective (and, often, solid-on-solid) models, 
in which the interface is still microscopic---represented by a stochastic field---while the structural details of the bulk phases are neglected.  

A simple example of such an effective model is a \emph{gradient field}. 
To define this system, we consider a finite subset~$\Lambda$ of the $d$-dimensional hypercubic lattice~$\Z^d$ and, 
at each site of~$\Lambda$ and its external boundary~$\partial\Lambda$, 
we consider the real-valued variable~$\phi_x$ representing the height of the interface at~$x$. 
The Hamiltonian is then given by
\begin{equation}
\label{1.1} 
H_\Lambda(\phi)=\sum_{\begin{subarray}{c}
\langle x,y\rangle\\ x\in\Lambda, y\in\Lambda\cup\partial\Lambda
\end{subarray}} V(\phi_y-\phi_x),
\end{equation} 
where the sum is over unordered nearest-neighbor pairs $\langle x,y\rangle$. 
A standard example is the quadratic potential~$V(\eta)=\frac12\kappa\eta^2$ with~$\kappa>0$; 
in general~$V$ is assumed to be a smooth, even function with 
a sufficient (say, quadratic) growth at infinity. 
The Gibbs measure takes the usual form
\begin{equation}
P_\Lambda(\textd\phi)=Z^{-1}\texte^{-\beta H_\Lambda(\phi)}\textd\phi,
\end{equation} 
where~$\textd\phi$ is the $|\Lambda|$-dimensional Lebesgue measure (the boundary values of~$\phi$ remain fixed and implicit in
the notation), $\beta>0$ is the inverse temperature and~$Z$ is a normalization constant.

A natural question to ask is what are the possible 
limits of the Gibbs measures $P_\Lambda(\textd\phi)$ as~$\Lambda\uparrow\Z^d$.
Unfortunately, in dimensions~$d=1,2$, the fields~$(\phi_x)_{x\in\Lambda}$ are 
very ``rough'' 
no matter how tempered the boundary conditions are assumed to be. 
As a consequence, the family of measures~$(P_\Lambda)_{\Lambda\subset\Z^d}$ is not tight 
and no meaningful object is obtained by taking the  
limit~$\Lambda\uparrow\Z^d$---i.e., the interface is \emph{delocalized}.
On the other hand, in dimensions~$d\ge3$ the fields are sufficiently smooth to permit a non-trivial thermodynamic limit---the interface is \emph{localized}. 
These facts are established by combinations of Brascamp-Lieb inequality techniques and/or random walk representation (see, e.g.,~\cite{Funaki})  which, unfortunately, 
apply only for convex potentials with uniformly positive curvature. 
Thus,
somewhat surprisingly, even for $V(\eta)=\eta^4$ the problem of localization in high-dimension is still 
open~\cite[Open Problem~1]{Velenik}. 

As it turns out, the thermodynamic limit of the measures~$P_\Lambda$ is significantly less singular 
once we restrict attention to the gradient variables
$\eta=(\eta_b)$.
These are defined by~$\eta_b=\phi_y-\phi_x$ where~$b$ is the nearest-neighbor edge~$(x,y)$ oriented 
in one of the positive lattice directions.
Indeed, the~$\eta$-marginal of~$P_\Lambda(\textd\phi)$ always has at least one (weak) limit 
``point''   as~$\Lambda\to\Z^d$. 
The limit measures satisfy a natural DLR condition and are therefore called \emph{gradient Gibbs measures}.
(Precise definitions will be stated below or can be found in~\cite{Funaki,Sheffield}.)  
One non-standard aspect of the gradient variables is that they have to obey a host of constraints. 
Namely, 
\begin{equation}
\label{E:loop}
\eta_{b_1}+\eta_{b_2}=\eta_{b_3}+\eta_{b_4}
\end{equation}
holds for each lattice plaquette~$(b_1,b_2,b_3,b_4)$, 
where the edges~$b_j$ are listed counterclockwise and are assumed to be positively oriented. 
These constraints will be implemented at the level of \emph{a priori} measure, see Sect.~\ref{S:model}.

It would be natural to expect that the character (and number) of gradient Gibbs measures depends
sensitively on the potential~$V$. 
However,  this is not the case for the class of
uniformly strictly-convex potentials  (i.e., the~$V$'s such that~$V''(\eta)\ge c_->0$ for all~$\eta$).
Indeed, Funaki and Spohn \cite{Funaki-Spohn} 
showed that, in these cases, the translation-invariant, ergodic, gradient Gibbs measures are completely characterized 
by the \emph{tilt} of the underlying interface. 
Here
the tilt is a vector~$u\in\R^d$ such that
\begin{equation}
\label{E:tilt}
E\eta_b=u\cdot b
\end{equation} 
for every edge~$b$---which we regard as a vector in~$\R^d$.
Furthermore,
the correspondence is one-to-one, i.e., for each tilt there exists precisely one gradient Gibbs measure with this tilt. 
Alternative proofs permitting extensions to discrete gradient models have appeared in 
Sheffield's thesis~\cite{Sheffield}.  

It is natural to expect that a serious violation of the strict-convexity assumption on~$V$ may invalidate the above results.
Actually, an example of a gradient model with multiple gradient Gibbs states of the same tilt has recently been 
presented~\cite{Sheffield};  unfortunately, 
the example is not of the type considered above because of the lack of translation invariance and its reliance on the discreteness of the fields. 
The goal of this paper is to point out a general mechanism by which the model \eqref{1.1} with a
sufficiently non-convex potential~$V$ fails the conclusions of Funaki-Spohn's~theorems.

\subsection{Potentials of interest}
\label{S:models}\noindent
The mechanism driving our example will be the occurrence of a structural surface phase transition. To motivate the forthcoming considerations,
let us recall that phase transitions typically arise via one of two mechanisms: either due to the breakdown of an internal symmetry, or via an
abrupt turnover between energetically and entropically favored states. The standard examples of systems with these kinds of phase transitions are
the Ising model and the $q$-state Potts model with a sufficiently large~$q$, respectively.  In the former, at sufficiently low temperatures, there is
a spontaneous breaking of the symmetry between the plus and minus spin states; in the latter, there is a first-order transition at intermediate
temperatures between $q$~ordered, low-temperature states and a disordered, high-temperature state.

Our goal is to come up with a potential~$V$  that
would mimic one of the above situations. In the present context an analogue of
the Ising model appears to be  
a \emph{double-well potential} of the form,~e.g.,
\begin{equation}
\label{double-well}
V(\eta)=\kappa(\eta^2-\eta_\star^2)^2.
\end{equation}
Unfortunately, due to the underlying plaquette constraints \eqref{E:loop},  the symmetry between the wells cannot be completely broken and, even at the level of ground states, 
the system appears to be disordered. On~$\Z^2$ this can be demonstrated explicitly by making a link to the \emph{ice model}, which is a special case of the six vertex model~\cite{Baxter}. 
A similar equivalence has been used~\cite{vB} to study a roughening transition in an SOS~interface.

\begin{figure}[t]
\bigskip
\centerline{\epsfig{figure=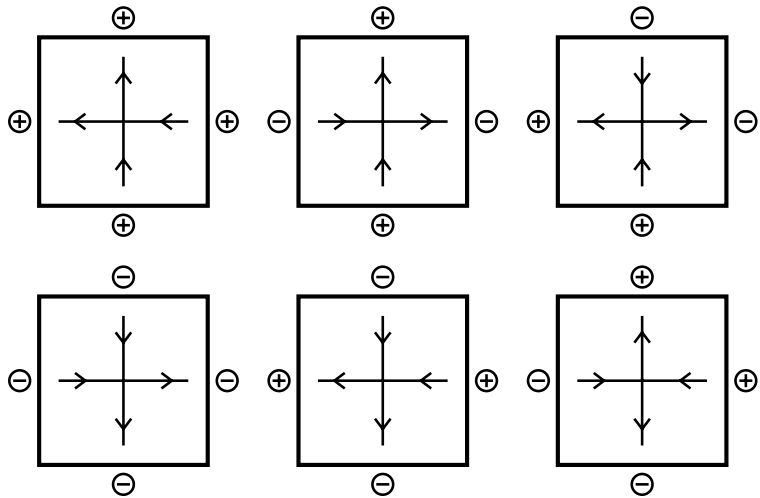, width=0.9\textwidth}}
\caption{The six plaquette configurations of minimal energy for the potential~\eqref{double-well} on~$\Z^2$ and their equivalent ice model configurations at the corresponding vertex on the dual grid. The sign marks represent the signs of~$\eta_b$ along the side of the plaquette $(b_1,b_2,b_3,b_4)$---with horizontal bonds $b_1,b_3$ oriented to the right and vertical bonds $b_2,b_4$ oriented upwards. The unit flow represented by the arrows runs upwards (downwards) through horizontal bonds with positive (negative) sign, and to the left (right) through vertical bonds with positive (negative) sign. The loop condition \eqref{E:loop} makes the flow conserved (i.e., no sources or sinks).
}
\label{fig1}
\end{figure}

To see how the equivalence works exactly, note that the ground states of the system \eqref{double-well} are such that all~$\eta$'s equal~$\pm\eta_\star$. Let us associate a unit flow with each \emph{dual} bond whose sign is determined by the value of~$\eta_b$ for its direct counterpart~$b$.  
The plaquette constraint \eqref{E:loop} then translates into a \emph{no-source-no-sink} condition for this flow. If we mark the flow by arrows, the dual bonds at each plaquette are constrained to one of six zero-flux arrangements of the six vertex model; cf~Fig.~\ref{fig1} and its caption.
The weights of all zero-flux arrangements are equal; we thus have the special case corresponding to the ice model. 
The ice model can be ``exactly  solved''~\cite{Baxter}: The ground states have a non-vanishing residual entropy~\cite{Lieb} and are disordered with infinite correlation length~\cite[Sect. 8.10.III]{Baxter}. However, it is not clear how much of this picture survives to positive temperatures.

\begin{figure}[t]
\bigskip
\centerline{\epsfig{figure=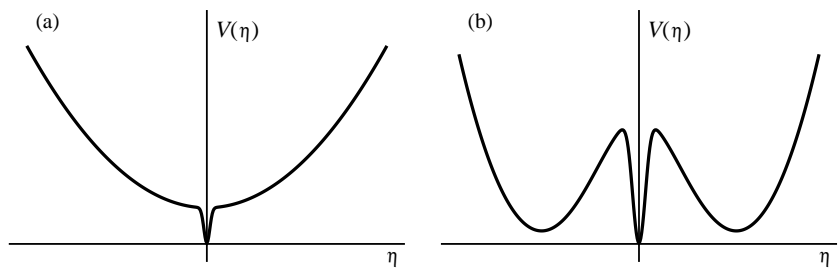, width=\textwidth}}
\caption{Two canonical examples of potentials that will lead to a structural surface phase transition. The picture labeled (a) is obtained by superimposing---in the sense of \eqref{V-double}---two symmetric wells of (significantly) different widths. Part~(b) of the figure represents the triple-well potential as defined in \eqref{V-triple}. For the application of our technique of proof, it only matters that the widths of the wells are sufficiently different.
}
\label{fig2}
\end{figure}

The previous discussion shows that it will be probably quite hard to realize a symmetry-breaking transition in the context of the gradient model~\eqref{1.1}. It is the order-disorder mechanism for phase transitions that seems considerably more promising.
There are two canonical examples of interest: a potential with \emph{two centered wells} and a \emph{triple-well potential}; see Fig.~\ref{fig2}. Both of these lead to a gradient model which features a phase transition, at some intermediate temperature, from states with the~$\eta$'s lying (mostly) within the thinner well to states whose~$\eta$'s fluctuate on the scale of the thicker well(s).

Our techniques apply equally to these---as well as other similar---cases provided the widths of the wells are sufficiently distinct. Notwithstanding, the analysis becomes significantly cleaner if we abandon temperature as our principal parameter (e.g., we set~$\beta=1$) and consider potentials~$V$ that are simply \emph{defined} by
\begin{equation}
\label{V-double}
\texte^{-V(\eta)}=p\,\texte^{-\kappaO\eta^2/2}+(1-p)\,\texte^{-\kappaD\eta^2/2}.
\end{equation}
Here~$\kappaO$ and~$\kappaD$ are positive numbers and~$p$ is a parameter taking values in~$[0,1]$. For appropriate values of the constants,~$V$ defined this way will have a graph as in Fig.~\ref{fig2}(a). To get the graph in part~(b), we would need to consider~$V$'s of the form
\begin{equation}
\label{V-triple}
\texte^{-V(\eta)}=p\,\texte^{-\kappaO\eta^2/2}+\frac{1-p}2\,\texte^{-\kappaD(\eta-\eta_\star)^2/2}
+\frac{1-p}2\,\texte^{-\kappaD(\eta+\eta_\star)^2/2},
\end{equation}
where~$\pm\eta_\star$ are the (approximate) locations of the off-center wells.

The idea underlying the expressions \eqref{V-double} and \eqref{V-triple} is similar to that of the Fortuin-Kasteleyn representation of the Potts model~\cite{FK}.
In the context of conti\-nuous-spin models similar to ours, such representation has fruitfully been used by Zahradn{\'\i}k~\cite{Z}. 
Focusing on \eqref{V-double}, we can interpret the terms on the right-hand side of 
\eqref{V-double} as two distinct states of each bond. (We will soon exploit this interpretation in detail.) The indexing of the coupling constants suggests the names: ``O'' for \emph{ordered} and ``D'' for \emph{disordered}. 

It is clear that the extreme values of~$p$ (near zero or near one) will be dominated by one type of bonds; what we intend to show is that, for $\kappaO$ and~$\kappaD$ sufficiently distinct from each other, the transition between the ``ordered'' and ``disordered'' phases is (strongly) first order. Similar conclusions and proofs---albeit more complicated---apply also to the potential~\eqref{V-triple}. However, for clarity of exposition, we will focus on the potential \eqref{V-double} for the rest of the paper
(see, however, Sect.~\ref{sec2.5}). 
In addition, we will also restrict ourselves to two dimensions, even though the majority of our results are valid for all~$d\ge2$.

\section{Main results}
\label{S:model}
\subsection{Gradient Gibbs measures}
We commence with 
a precise definition of our model. Most of the work in this paper will be confined to the lattice torus~$\T_L$ of~$L\times L$ sites in~$\Z^2$, so we will start with this particular geometry. Choosing the natural positive direction for each lattice axis, let~$\B_L$ denote the corresponding set of positively oriented edges in~$\T_L$. Given a configuration~$(\phi_x)_{x\in\T_L}$, we introduce the gradient  field~$\eta=\nabla\phi$ by assigning the variable $\eta_b=\phi_y-\phi_x$ to each~$b=(x,y)\in \B_L$. 
The product Lebesgue measure~$\prod_{x\ne0}\textd\phi_x$ induces a ($\sigma$-finite) measure~$\nu_L$ 
on the space~$\R^{\B_L}$ via
\begin{equation}
\label{E:nuL}
\nu_L\bigl(\AA)=\int\biggl(\,\prod_{x\in\T_L\smallsetminus\{0\}}\textd\phi_x\biggr)\,\delta(\textd\phi_0)\,
\1_{\{\nabla\phi\in\AA\}},
\end{equation}
where~$\delta$ denotes the Dirac point-mass at zero.

We interpret  the measure~$\nu_L$ as an \emph{a priori} measure on \emph{gradient} configurations~$\eta\in \R^{\B_L}$. 
Since the~$\eta$'s arise as the gradients of the~$\phi$'s it is easy to check that~$\nu_L$ is entirely supported on 
the linear subspace $\XX_L\subset \R^{\B_L}$ of configurations determined by the condition that the sum of signed~$\eta$'s---with a positive or negative sign depending on whether the edge is traversed in the positive or negative direction, respectively---vanishes  around each closed circuit on~$\T_L$.
(Note that, in addition to \eqref{E:loop}, the condition includes also loops that wrap around the torus.)
We will refer to such configurations as \emph{curl-free}.

Next we will define gradient Gibbs measures on~$\T_L$. For later convenience we will proceed in some more generality than presently needed: Let~$(V_b)_{b\in\B_L}$ be a collection of measurable functions $V_b\colon \R\to [0,\infty)$ and consider the partition function
\begin{equation}
Z_{L,(V_b)}=\int_{\R^{|\B_L|}}\exp\Bigl\{-\sum_{b\in\B_L} V_b(\eta_b)\Bigr\}\,\nu_L(\textd\eta). 
\end{equation}
Clearly, $Z_{L,(V_b)}>0$ and, under the condition that~$\eta\mapsto\texte^{-V_b(\eta)}$ is integrable 
with respect to the Lebesgue measure on~$\R$, also~$Z_{L,(V_b)}<\infty$. 
We may then define~$P_{L,(V_b)}$ to be the probability measure on~$\R^{|\B_L|}$ given by
\begin{equation}
\label{E:Pmeasure}
P_{L,(V_b)}(\textd\eta)=\frac{1}{Z_{L,(V_b)}}\exp\Bigl\{-\sum_{b\in\B_L} V_b(\eta_b)\Bigr\}\,\nu_L(\textd\eta).
\end{equation}
This is the \emph{gradient Gibbs measure} on~$\T_L$ corresponding to the potentials~$(V_b)$. In the situations when 
$V_b=V$ for all~$b$---which is the principal case of interest in this paper---we will denote the corresponding gradient Gibbs measure on~$\T_L$ by~$P_{L,V}$.

It is not surprising that~$P_{L,(V_b)}$ obeys appropriate DLR equations with respect to all connected $\Lambda\subset\T_L$  containing no topologically non-trivial circuit.
Explicitly, if~$\eta_{\Lambda^\cc}$ in~$\Lambda^\cc$ is a curl-free boundary condition, then 
the conditional law of~$\eta_\Lambda$ given~$\eta_{\Lambda^\cc}$ is
\begin{equation}
\label{specifikace}
P_{L,(V_b)}(\textd\eta_\Lambda|\eta_{\Lambda^\cc})
=\frac1{Z_\Lambda(\eta_{\Lambda^\cc})}\exp\Bigl\{-\sum_{b\in\Lambda}V_b(\eta_b)\Bigr\}\,
\nu_\Lambda(\textd\eta_\Lambda|\eta_{\Lambda^\cc}).
\end{equation}
Here $P_{L,(V_b)}(\textd\eta_\Lambda|\eta_{\Lambda^\cc})$ is the conditional probability 
with respect to the (tail)~$\sigma$-algebra~$\scrT_\Lambda$ 
generated 
by the fields on $\Lambda^\cc$, $Z_\Lambda(\eta_{\Lambda^\cc})$ is the partition function in~$\Lambda$, and~$\nu_\Lambda(\textd\eta_\Lambda|\eta_{\Lambda^\cc})$ is the \emph{a priori} measure induced by $\nu_L$ on~$\eta_\Lambda$ given the boundary condition~$\eta_{\Lambda^\cc}$. 

As usual, this property remains valid even in the thermodynamic limit.
We thus say that a measure on~$\mu$ is an \emph{infinite-volume gradient Gibbs measure} if it satisfies the DLR equations with respect to the specification \eqref{specifikace} in any finite set~$\Lambda\subset\Z^2$. (As is easy to check---e.g., by reinterpreting the~$\eta$'s back in terms of the~$\phi$'s---$\nu_\Lambda(\textd\eta_\Lambda|\eta_{\Lambda^\cc})$ is independent of the values of~$\eta_{\Lambda^\cc}$ outside any circuit winding around~$\Lambda$, and so it is immaterial that it originated from a measure on torus.)

\smallskip
An important aspect of our derivations will be the fact that our potential~$V$ takes the specific form \eqref{V-double}, which can 
be concisely written as
\begin{equation}
\label{V-rep}
\texte^{-V(\eta)}=\int\varrho(\textd\kappa)\,\texte^{-\frac12\kappa\eta^2},
\end{equation}
where~$\varrho$ is the probability measure 
$\varrho=p\delta_{\kappaO}+(1-p)\delta_{\kappaD}$. It follows that the Gibbs measure~$P_{L,V}$ can be regarded as the projection of the 
\emph{extended gradient Gibbs measure},
\begin{equation}
\label{E:Qmeasure}
Q_L(\textd\eta,\textd\kappa)
=\frac{1}{Z_{L,V}}\exp\Bigl\{-\frac12\sum_{b\in\B_L} \kappa_b\eta_b^2\Bigr\}\,\nu_L(\textd\eta)\varrho_L(\textd\kappa),
\end{equation}
to the $\sigma$-algebra generated by the~$\eta$'s. 
Here~$\varrho_L$ is the product of measures~$\varrho$, one for each bond in~$\B_L$. 
As is easy to check, conditioning on $(\eta_b, \kappa_b)_{b\in\Lambda^\cc}$ yields the corresponding  extension $Q_\Lambda(\textd\eta_\Lambda\textd\kappa_\Lambda|\eta_{\Lambda^\cc})$ of the finite-volume specification \eqref{specifikace}---the result is independent of the~$\kappa$'s outside~$\Lambda$ because, once~$\eta_{\Lambda^\cc}$ is fixed, 
these have no effect on the configurations~in~$\Lambda$.

The main point of introducing the extended measure is that, if conditioned on the~$\kappa$'s, the variables~$\eta_b$ are distributed as gradients of a Gaussian field---albeit with a non-translation invariant covariance matrix. As we will see, the phase transition proved in this paper is manifested by a jump-discontinuity in the density of bonds with~$\kappa_b=\kappaO$ 
which at the level of~$\eta$-marginal results in
a jump in the characteristic scale of the fluctuations.

\begin{myremark}
Notably, the extended measure~$Q_L$ plays the same role for~$P_{L,V}$ as the so called Edwards-Sokal coupling measure~\cite{ES} does for the Potts model. 
Similarly as for the Edwards-Sokal measures~\cite{Grimmett, BBCK}, there is a one-to-one correspondence between the infinite-volume measures on~$\eta$'s and the corresponding 
infinite-volume extended gradient Gibbs
measures on~$(\eta,\kappa)$'s. Explicitly, if~$\mu$ is an infinite-volume gradient Gibbs measure for potential~$V$, then~$\tilde\mu$, defined by (extending the consistent family of measures of the~form)
\begin{equation}
\label{kappa-ext}
\tilde\mu\bigl((\eta_b,\kappa_b)_{b\in\Lambda}\in\AA\times\BB\bigr)
=\int_\BB\varrho_\Lambda(\textd\kappa)\,
E_\mu\biggl(\,\1_\AA\prod_{b\in\Lambda}\texte^{-\frac12\kappa_b\eta_b^2+V(\eta_b)}\biggr),
\end{equation}
is a Gibbs measure with respect to the extended specifications $Q_\Lambda(\cdot|\eta_{\Lambda^\cc})$. For the situations with only a few distinct values of~$\kappa_b$, it may be of independent interest to study the properties of the $\kappa$-marginal of the extended measure, e.g., using the techniques of percolation theory. However, apart from some remarks in Sect.~\ref{sec-diskuse}, we will not pursue these matters in the present paper.
\end{myremark}

\subsection{Phase coexistence of gradient measures}
\label{sec2.2}\noindent
Now we are ready to state our main results. Throughout we will consider the potentials~$V$ of the form \eqref{V-double} with~$\kappaO\gg\kappaD$. As a moment's thought reveals, the model is invariant under the transformation
\begin{equation}
\kappaO\to\kappaO\theta^2,\quad\kappaD\to\kappaD\theta^2,\quad\eta_b\to\eta_b/\theta
\end{equation}
for any fixed~$\theta\ne0$.
In particular, without loss of generality, one could assume from the beginning that~$\kappaO\kappaD=1$
and regard~$\ffrac\kappaO\kappaD$ as the sole parameter of the model. 
However, we prefer to treat the two terms in \eqref{V-double} on an equal footing, and so we will keep the coupling strengths independent.

Given a shift-ergodic gradient Gibbs measure, recall that its tilt is the vector~$u$ such that \eqref{E:tilt} holds for each bond. The principal result of the present paper is the following theorem:

\begin{theorem}
\label{T:main}
For each~$\epsilon>0$ there exists a constant~$c=c(\epsilon)>0$ and, if
\begin{equation}
\label{kappa-cond}
\kappaO\ge c\kappaD,
\end{equation}
a number~$\pt\in(0,1)$ such that,  for interaction~$V$ with~$p=\pt$, there are
two distinct, infinite-volume, shift-ergodic gradient Gibbs measures~$\mu_\ord$ and~$\mu_\dis$ of zero tilt for which 
\begin{equation}
\label{E:ener-bd1}
\mu_\ord\Bigl(\,|\eta_b|\ge\frac\lambda{\sqrt{\kappaO}}\Bigr)\le \epsilon+\frac1{4\lambda^2},
\qquad\forall\lambda>0,
\end{equation}
and
\begin{equation}
\label{E:ener-bd2}
\mu_\dis\Bigl(\,|\eta_b|\le\frac\lambda{\sqrt{\kappaD}}\Bigr)\le \epsilon+c_1\lambda^{\ffrac14},
\qquad\forall\lambda>0.
\end{equation}
Here~$c_1$ is
a constant of order unity.
\end{theorem}

\begin{myremark}
An inspection of the proof actually reveals that the above bounds are valid for any~$\epsilon$ satisfying~$\epsilon\ge c_2(\ffrac\kappaD\kappaO)^{1/8}$, where~$c_2$ is a constant of order unity.
\end{myremark}

As already alluded to, this result is a consequence of the fact that the density of ordered bonds, i.e., those with~$\kappa_b=\kappaO$, undergoes a jump at~$p=\pt$. On the torus, we can make the following asymptotic statements:

\begin{theorem}
\label{T:torus}
Let~$R^\ord_L$ denote the fraction of ordered bonds on~$\T_L$, i.e.,
\begin{equation}
\label{Rord}
R^\ord_L=\frac1{|\B_L|}\sum_{b\in\B_L}\1_{\{\kappa_b=\kappaO\}}.
\end{equation}
For each~$\epsilon>0$ there exists~$c=c(\epsilon)>0$ such that the following holds:
Under the condition~\eqref{kappa-cond}, and for~$\pt$ as in 
Theorem~\ref{T:main},
\begin{equation}
\label{below-pt}
\lim_{L\to\infty}Q_L(R^\ord_L\le\epsilon)=1,\qquad p<\pt
\end{equation}
and
\begin{equation}
\label{above-pt}
\lim_{L\to\infty}Q_L(R^\ord_L\ge1-\epsilon)=1,\qquad p>\pt.
\end{equation}
\end{theorem}

The present setting actually permits us to determine the value of~$\pt$ via a duality argument. This is the only result in this paper which is intrinsically two-dimensional
(and intrinsically tied to the form \eqref{V-double} of~$V$). All 
other conclusions can be extended to~$d\ge2$ and 
to more general potentials.

\begin{theorem}
\label{T:dual}
Let $d=2$.
If\/ $\ffrac\kappaO\kappaD\gg1$, then~$\pt$ is given by
\begin{equation}
\label{pt-eq}
\frac{\pt}{1-\pt}=\Bigl(\frac{\kappaD}{\kappaO}\Bigr)^{\ffrac14}.
\end{equation}
\end{theorem}

Theorem~\ref{T:torus} is proved in Sect.~\ref{sec4.2}, Theorem~\ref{T:main} is proved in Sect.~\ref{sec4.3} and Theorem~\ref{T:dual} is proved in Sect.~\ref{sec5.3}.

\subsection{Discussion}
\label{sec-diskuse}\noindent
The phase transition described in the above theorems can be interpreted in several ways. First, in terms of the extended gradient Gibbs measures on torus, it clearly corresponds to a transition between a state with nearly all bonds ordered ($\kappa_b=\kappaO$) to a state with nearly all bonds disordered ($\kappa_b=\kappaD$). Second, looking back at the inequalities \twoeqref{E:ener-bd1}{E:ener-bd2}, most of the~$\eta$'s will be of order at most~$1/\sqrt{\kappaO}$ in the ordered state while most of them will be of order at least~$1/\sqrt{\kappaD}$ in the disordered state. Hence, the corresponding (effective) interface is significantly rougher at~$p<\pt$ than it is at~$p>\pt$ (both phases are rough according to the standard definition of this term) and we may thus interpret the above as a kind of \emph{first-order roughening} transition that the interface undergoes at~$\pt$. Finally, since the gradient fields in the two states fluctuate on different characteristic scales, the entropy (and hence the energy) associated with these states is different; we can thus view this as a standard energy-entropy transition. (By the energy we mean the expectation of~$V(\eta_b)$; notably, the expectation of~$\kappa_b\eta_b^2$ is the same in both measures; cf~\eqref{E:var-bd}.) 
Energy-entropy transitions for spin models have been studied in~\cite{DS,KS,KS-proceedings} and, quite recently, in~\cite{vES}.

Next let us turn our attention to the conclusions of Theorem~\ref{T:torus}. We actually believe that the dichotomy \twoeqref{below-pt}{above-pt} applies (in the sense of almost-sure limit of~$R^\ord_L$ as~$L\to\infty$) to all
translation-invariant extended
gradient Gibbs states with zero tilt. The reason is that, conditional on the~$\kappa$'s, the gradient fields are Gaussian with uniformly positive stiffness. We rest assured that the techniques of~\cite{Funaki-Spohn} and~\cite{Sheffield} can be used to prove that the gradient Gibbs measure with zero tilt is unique for almost every
configuration of the~$\kappa$'s; so the only reason for multiplicity of gradient Gibbs measures with zero tilt is a phase transition in the~$\kappa$-marginal. However, a detailed write-up of this argument would require developing the  precise---and
somewhat subtle---correspondence between the gradient Gibbs measures of a given tilt and the minimizers of the Gibbs variational principle (which we have, in full detail, only for convex periodic potentials~\cite{Sheffield}). Thus, to keep the paper at manageable length, we limit ourselves to a weaker result.

The fact that the transition occurs at~$\pt$ satisfying \eqref{pt-eq} is a consequence of a \emph{duality} between the $\kappa$-marginals at~$p$ and~$1-p$. 
More generally, the duality links the marginal law of the configuration~$(\kappa_b)$ with the law of~$(1/\kappa_b)$; see Theorem~\ref{T:dualita} and Remark~\ref{remark:mira-dualita}. 
[At the level of gradient fields, the duality provides only a vague link between the flow of the weighted gradients~$(\sqrt{\kappa_b}\eta_b)$ along a given curve and its flux through this curve. Unfortunately, this link does not seem to be particularly useful.] The point~$p=\pt$ is self-dual which makes it the most natural candidate for a transition point. It is interesting to ponder about what happens when~$\ffrac\kappaO\kappaD$ decreases to one. Presumably, the first-order transition (for states at zero tilt) disappears before~$\ffrac\kappaO\kappaD$ reaches one and is replaced by some sort of critical behavior. 
Here the first problem to tackle is to establish the \emph{absence} of first-order phase transition for small~$\ffrac\kappaO\kappaD-1$. Via a standard duality argument (see~\cite{Chayes-Shtengel}) this would yield a power-law lower bound for bond connectivities at~$\pt$.

Another interesting problem is to determine what happens with measures of non-zero tilt. We expect that, at least for moderate values of the tilt~$u$, the first-order transition persists but shifts to lower values of~$p$. Thus, one could envision a whole phase diagram in the $p$-$u$ plane. Unfortunately, we are unable to make any statements of this kind because the standard ways to induce a tilt on the torus (cf~\cite{Funaki-Spohn}) lead to measures that are not reflection positive.

\subsection{Outline of the proof}
We proceed by an outline of the principal steps of the proof to which the remainder of this paper is devoted. The arguments are close in spirit to those in~\cite{DS,KS,KS-proceedings}; the differences arise from the subtleties in the setup due to the gradient nature of the fields. 

The main line of reasoning is basically thermodynamical: Consider the $\kappa$-marginal of the extended torus state~$Q_L$ which we will regard as a measure on configurations of ordered and disordered bonds. Let~$\chi(p)$ denote (the $L\to\infty$ limit of) the expected fraction of ordered bonds in the torus state at parameter~$p$. Clearly~$\chi(p)$ increases from zero to one as~$p$ sweeps through~$[0,1]$. The principal observation is that, under the assumption~$\ffrac\kappaO\kappaD\gg1$, the quantity~$\chi(1-\chi)$ is small, uniformly in~$p$. Hence, $p\mapsto\chi(p)$ must undergo a jump from values near zero to values near one at some~$\pt\in(0,1)$. By usual weak-limiting arguments we construct two distinct gradient Gibbs measures at~$\pt$, one with high density of ordered bonds and the other with high density of disordered bonds.

The crux of the matter is thus to justify the uniform smallness of~$\chi(1-\chi)$. This will be a consequence of the fact that the simultaneous occurrence of ordered and disordered bonds at any two given locations is (uniformly) unlikely. For instance, let us estimate the probability that a particular plaquette has two ordered bonds emanating out of one corner and two disordered bonds emanating out of the other. Here the technique of chessboard estimates~\cite{FL,FILS1,FILS2} allows us to disseminate this pattern all over the torus via successive reflections (cf Theorem~\ref{T:chess} in Sect.~\ref{S:chessboard}). This bounds the quantity of interest by the~$1/L^2$-power of the probability that every other horizontal (and vertical) line is entirely ordered and the remaining lines are disordered. The resulting ``spin-wave calculation''---i.e., diagonalization of a period-2 covariance matrix in the Fourier basis and taking its determinant---is performed (for all needed patterns) in Sect.~\ref{S:spin-wave}.

Once the occurrence of a ``bad pattern'' is estimated by means of various spin-wave free energies, we need to prove that these ``bad-pattern'' spin-wave free energies are always worse off than those of the homogeneous patterns (i.e., all ordered or all disordered)---this is the content of Theorem~\ref{T:min}. Then we run a standard Peierls' contour estimate whereby the smallness of~$\chi(1-\chi)$ follows. Extracting two distinct, infinite-volume, ergodic gradient Gibbs states~$\mu_\ord$ and~$\mu_\dis$ at~$p=\pt$, it remains to show that these are both of zero tilt. Here we use the fact that, conditional on the~$\kappa$'s, the torus measure is symmetric Gaussian with uniformly positive stiffness. Hence, we can use standard Gaussian inequalities to show exponential tightness of the tilt, uniformly in the~$\kappa$'s; cf Lemma~\ref{L:tilt}. Duality calculations (see Sect.~\ref{S:duality}) then yield~$p=\pt$.

\subsection{Generalizations}
\label{sec2.5}\noindent
Our proof of phase coexistence applies to any potential of the form 
shown 
in Fig.~\ref{fig2}---even if we return to parametrization by~$\beta$. The difference with respect to the present setup is that in the general case we would have to approximate the potentials by a quadratic well at each local minimum and, before performing the requisite Gaussian calculations, estimate the resulting~errors.

Here is a sketch of the main ideas: We fix a scale~$\Delta$ and regard~$\eta_b$ to be in a well if it is within~$\Delta$ of the corresponding local minimum. Then the requisite quadratic approximation of~$\beta$-times energy is good up to errors of order~$\beta\Delta^3$. The rest of the potential ``landscape'' lies at energies of at least order~$\Delta^2$ and so it will be only ``rarely visited'' by the~$\eta$'s provided that~$\beta\Delta^2\gg1$. On the other hand, the same condition ensures that the spin-wave integrals are essentially not influenced by the restriction that~$\eta_b$ be within~$\Delta$ of the local minimum. Thus, to make all approximations work we need that
\begin{equation}
\beta \Delta^3\ll1\ll\beta\Delta^2
\end{equation}
which is achieved for~$\beta\gg1$ by, e.g.,~$\Delta=\beta^{-\frac5{12}}$.
This approach has recently been used to prove phase transitions in classical~\cite{BCN,BCKiv} as well as quantum~\cite{BCStarr} systems with highly degenerate ground states. We refer the reader to these references for further details.

A somewhat more delicate issue is the proof that both coexisting states are of zero tilt. Here the existing techniques require that we have some sort of uniform convexity. This more or less forces us to use the~$V$'s of the form
\begin{equation}
V(\eta)=-\log\biggl(\,\sum_j\texte^{-V_j(\eta)}\biggr),
\end{equation}
where the~$V_j$'s are uniformly convex functions. Clearly, our choice~\eqref{V-double} is the simplest potential of this type; the question is how general the potentials obtained this way can be. We hope to return to this question in a future publication.

\section{Spin-wave calculations}
\label{S:spin-wave}\noindent
As was just mentioned, the core of our proofs are estimates of the spin-wave free energy for various regular patterns of ordered and disordered bonds on the torus. These estimates are rather technical and so we prefer to clear them out of the way before we get to the main line of the proof. The readers wishing to follow the proof in linear order may consider skipping this section and returning to it only while reading the arguments in Sect.~\ref{sec4.2}. Throughout this and the forthcoming sections we assume that~$L$ is an even integer.

\subsection{Constrained partition functions}
\label{S:SWFE}
\begin{figure}[t]
\bigskip
\centerline{\epsfig{figure=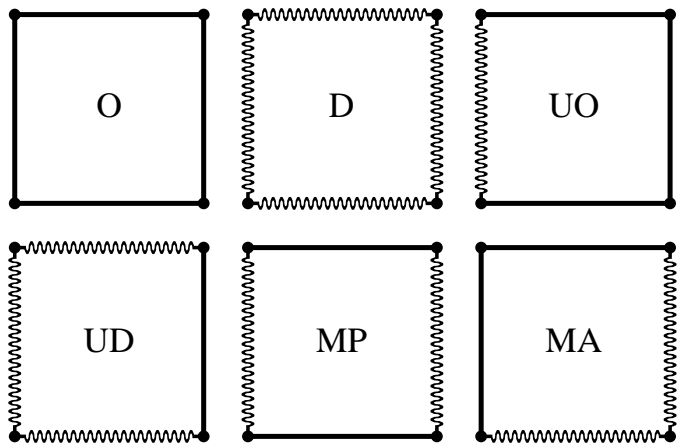, width=3.6in}}
\caption{Six possible arrangements of ``ordered'' and ``disordered'' bonds around a lattice plaquette. Here the ``ordered'' bonds are represented by solid lines and the ``disordered'' bonds by wavy lines. Each inhomogeneous pattern admits other rotations which are not depicted. The acronyms stand for: (top row) Ordered, Disordered and U-shape Ordered and (bottom row) U-shape Disordered, Mixed Periodic and Mixed Aperiodic, respectively.
}
\label{fig3}
\end{figure}

\noindent
We will consider six partition functions~$\ZO_L$, $\ZD_L$, $\ZUO_L$, $\ZUD_L$, $\ZMP_L$ and $\ZMA_L$ on~$\T_L$ that 
correspond to six regular configurations each of which is obtained by reflecting one of six possible arrangements of ``ordered'' and ``disordered'' bonds around a lattice plaquette to the entire torus. These quantities will be the ``building blocks'' of our analysis in Sect.~\ref{sec:PT}. The six plaquette configurations are depicted in Fig.~\ref{fig3}.

We begin by considering the homogeneous configurations. Here~$\ZO_L$ is the partition function~$Z_{L,(V_b)}$ for all edges of the ``ordered'' type:
\begin{equation}
V_b(\eta)=-\log p+\frac12\kappaO\eta^2,
\qquad b\in\B_L.
\end{equation}
Similarly,~$\ZD_L$ is the quantity~$Z_{L,(V_b)}$ for
\begin{equation}
V_b(\eta)=-\log(1-p)+\frac12\kappaD\eta^2,
\qquad b\in\B_L,
\end{equation}
i.e., with all edges ``disordered.''

Next we will define the partition functions~$\ZUO_L$ and~$\ZUD_L$ which are obtained by reflecting a plaquette with three bonds of one type and the remaining bond of the other type. Let us split~$\B_L$ into the even~$\Beven_L$ and odd~$\Bodd_L$ horizontal and vertical edges---with the even edges on the lines of sites in the~$x$ direction with even~$y$ coordinates and lines of sites in~$y$ direction with even~$x$ coordinates. 
Similarly, we will also consider the decomposition of~$\B_L$ into the set of horizontal edges~$\Bhor_L$ and vertical edges~$\Bvert_L$.
Letting
\begin{equation}
V_b(\eta)=\begin{cases}
-\log p+\frac12\kappaO\eta^2,\qquad&\text{if }b\in\Bhor_L\cup
\Beven_L,
\\
-\log(1-p)+\frac12\kappaD\eta^2,\qquad&\text{otherwise},
\end{cases}
\end{equation}
the partition function~$\ZUO_L$ then corresponds to the quantity~$Z_{L,(V_b)}$. The partition function~$\ZUD_L$ is obtained similarly; with the roles of ``ordered'' and ``disordered'' interchanged. Note that, since we are working on a square torus, the orientation of the pattern we choose does not matter.

It remains to define the partition functions~$\ZMP_L$ and~$\ZMA_L$ corresponding to the patterns with two ``ordered'' and two ``disordered'' bonds. For the former, we simply take~$Z_{L,(V_b)}$ with the potential
\begin{equation}
V_b(\eta)=\begin{cases}
-\log p+\frac12\kappaO\eta^2,\qquad&\text{if }b\in\Bhor_L,
\\
-\log(1-p)+\frac12\kappaD\eta^2,\qquad&\text{if }b\in\Bvert_L.
\end{cases}
\end{equation}
Note that the two types of bonds are arranged in a ``mixed periodic'' pattern; hence the index~$\MP$. As to the quantity~$\ZMA_L$, here we will consider a ``mixed aperiodic'' pattern. Explicitly, we define
\begin{equation}
V_b(\eta)=\begin{cases}
-\log p+\frac12\kappaO\eta^2,\qquad&\text{if }b\in\Beven_L,
\\
-\log(1-p)+\frac12\kappaD\eta^2,\qquad&\text{if }b\in\Bodd_L.
\end{cases}
\end{equation}
The ``mixed aperiodic'' partition function~$\ZMA_L$ is the quantity~$Z_{L,(V_b)}$ for this choice of~$(V_b)$. Again, on a square torus it is immaterial for the values of $\ZMP_L$ and~$\ZMA_L$ which orientation of the initial plaquette we start with.

As usual, associated with these partition functions are the corresponding free energies. In finite volume, these quantities can be defined in all cases by the formula
\begin{equation}
\label{finite-FE}
F_{L,\alpha}(p)=-\frac1{L^2}\log\frac{Z_{L,\alpha}}{(2\pi)^{\frac12(L^2-1)}},
\qquad \alpha=\Oind,\Dind,\UO,\UD,\MP,\MA,
\end{equation}
where the factor $(2\pi)^{\frac12(L^2-1)}$ has been added for later convenience and where the~$p$-dependence arises via the corresponding formulas for~$V_b$ in each particular case.

\subsection{Limiting free energies}
The goal of this section is to compute the thermodynamic limit of the~$F_{L,\alpha}$'s. 
For homogeneous and isotropic configurations, an important role will be played by the momentum representation of the lattice Laplacian
$\widehat D(\bk)=|1-\texte^{\texti k_1}|^2+|1-\texte^{\texti k_2}|^2$
defined for all~$\bk=(k_1,k_2)$ in the corresponding Brillouin zone~$\bk\in[-\pi,\pi]\times[-\pi,\pi]$.
Using this quantity, the ``ordered'' free energy will be simply
\begin{equation}
\label{E:FS}
F_{\Oind}(p)=-2\log p+\frac12\int_{[-\pi,\pi]^2}\frac{\textd\bk}{(2\pi)^2}\log\bigl\{\kappaO\widehat D(\bk)\bigr\},
\end{equation}
while the disordered free energy boils down to
\begin{equation}
\label{E:FA}
F_{\Dind}(p)=-2\log(1-p)
+\frac12\int_{[-\pi,\pi]^2}\frac{\textd\bk}{(2\pi)^2}\log\bigl\{\kappaD\widehat D(\bk)\bigr\}.
\end{equation}
It is easy to check that, despite the logarithmic singularity at~$\bk=\bzero$, both integrals converge.
The bond pattern underlying the quantity~$\ZMP_L$ lacks rotation invariance and so a different propagator 
appears inside the momentum integral:
\begin{multline}
\label{E:FMP}
F_{\MP}(p)=-\log\bigl[p(1-p)\bigr]\\
+\,\frac12\int_{[-\pi,\pi]^2}\frac{\textd\bk}{(2\pi)^2}\log\bigl\{\kappaO|1-
\texte^{\texti k_1}|^2+\kappaD|1-\texte^{\texti k_2}|^2\bigr\}.
\end{multline}
Again, the integral converges as long as (at least) one of~$\kappaO$ and~$\kappaD$ is strictly positive.

The remaining partition functions come from configurations that lack translation invariance and are ``only'' periodic with period two. 
Consequently, the Fourier transform of the corresponding propagator is only block diagonal, 
with two or four different~$\bk$'s ``mixed'' inside each block.
In the $\UO$ cases we will get the function
\begin{equation}
\label{E:FUO}
F_{\UO}(p)=-\frac12\log\bigl[p^3(1-p)\bigr]
\,+\,\frac14\int_{[-\pi,\pi]^2}\frac{\textd\bk}{(2\pi)^2}\log\bigl\{\det\Pi_{\UO}(\bk)\bigr\},
\end{equation}
where~$\Pi_{\UO}(\bk)$ is the $2\times2$-matrix
\begin{equation}
\label{PiUO}
\Pi_{\UO}(\bk)=\left(\,
\begin{matrix}
\kappaO|a_-|^2+\frac12(\kappaO+\kappaD)|b_-|^2 &
\frac12(\kappaO-\kappaD)
|b_-|^2
\\*[2mm]
\frac12(\kappaO-\kappaD)
|b_-|^2
&
\kappaO|a_+|^2+\frac12(\kappaO+\kappaD)
|b_-|^2
\end{matrix}
\,\right)
\end{equation}
with $a_\pm$ and~$b_\pm$ defined by
\begin{equation}
a_\pm=1\pm \texte^{\texti k_1}
\quad\text{and}\quad  b_\pm=1\pm \texte^{\texti k_2}.
\end{equation}
The extra factor~$\ffrac12$---on top of the usual~$\ffrac12$---in front of the integral arises because $\det\Pi_\UO(\bk)$ combines the contributions of two Fourier models; namely~$\bk$ and~$\bk+\pi\hate_1$.
A calculation shows
\begin{equation}
\label{detUO-bound}
\det\Pi_\UO(\bk)\ge \kappaO^2|a_-|^2|a_+|^2+\kappaO\kappaD|b_-|^4,
\end{equation}
implying that the integral in~\eqref{E:FUO} converges.
The free energy~$F_{\UD}$ is obtained by interchanging the roles of~$\kappaO$ and~$\kappaD$ and of~$p$ and~$(1-p)$.

In the~MA-cases we will assume that~$\kappaO\ne\kappaD$---otherwise there is no distinction between any of the six patterns. The corresponding free energy is then given~by
\begin{multline}
\label{E:FMA}
F_{\MA}(p)=-\log\bigl[p(1-p)\bigr]
\\+\,\frac18\int_{[-\pi,\pi]^2}\frac{\textd\bk}{(2\pi)^2}\log\biggl\{\Bigl(\frac{\kappaO-\kappaD}2\Bigr)^4\det\Pi_{\MA}(\bk)\biggr\}.
\end{multline}
Here $\Pi_{\MA}(\bk)$ is the $4\times4$-matrix
\begin{equation}
\label{E:Pi}
\Pi_{\MA}(\bk)=\left(\,
\begin{matrix}
\llambda(|a_-|^2+|b_-|^2)\!\!\! & |b_-|^2 & |a_-|^2 & 0
\\*[1mm] 
|b_-|^2 & \!\!\!\llambda(|a_+|^2+|b_-|^2) \!\!\!& 0 & |a_+|^2 
\\*[1mm]
|a_-|^2 & 0 & \!\!\!\llambda(|a_-|^2+|b_+|^2)\!\!\! & |b_+|^2
\\*[1mm] 
0 & |a_+|^2 & |b_+|^2 & \!\!\!\llambda(|a_+|^2+|b_+|^2)
\end{matrix}
\,\right)
\end{equation} 
with the abbreviation
\begin{equation}
\llambda=\frac{\kappaO+\kappaD}{\kappaO-\kappaD}.
\end{equation}
Note that~$r>1$ in the cases of our interest.

Observe that~$\det\Pi_\MA(\bk)$ is a quadratic polynomial in~$r^2$, i.e.,~$\det\Pi_\MA=Ar^4+Br^2+C$. Moreover,~$\Pi_\MA(\bk)$ annihilates~$(1,-1,-1,1)$ when~$r=1$, and so~$r^2=1$ is a root of~$Ar^4+Br^2+C$. Hence $\det\Pi_\MA(\bk)=(r^2-1)(Ar^2-C)$, i.e.,
\begin{multline}
\label{E:detPi}
\qquad
\det\Pi_{\MA}(\bk)=(\llambda^2-1)\Bigl\{ -(|a_+|^2|a_-|^2-|b_+|^2|b_-|^2)^2
\\
+(|a_+|^2+|b_+|^2)(|a_-|^2+|b_+|^2)(|a_+|^2+|b_-|^2)(|a_-|^2+|b_-|^2)\llambda^2
\Bigr\}.
\qquad
\end{multline} 
Setting~$r=1$ inside the large braces yields
\begin{equation}
\label{detMA-bound}
\det\Pi_{\MA}(\bk)\ge
4(\llambda^2-1)|a_-|^2|a_+|^2|b_-|^2|b_+|^2,
\end{equation}
implying that the integral in \eqref{E:FMA} is well defined and finite.

\begin{myremark}
The fact that~$\Pi_\MA(\bk)$ has zero eigenvalue at~$r=1$ is not surprising. Indeed,~$r=1$ corresponds to~$\kappaD=0$ in which case a quarter of all sites in the~$\MA$-pattern get decoupled from the rest. This indicates that the partition function blows up (at least) as~$(r-1)^{-|\T_L|/4}$ as~$r\downarrow1$ implying that there should be a zero eigenvalue at~$r=1$ per each~$4\times4$-block~$\Pi_\MA(\bk)$.
\end{myremark}

A formal connection between the quantities in \eqref{finite-FE} and those in \twoeqref{E:FS}{E:FMA} is guaranteed by the following result:

\begin{theorem}
\label{T:FE}
For all~$\alpha=\Oind,\Dind,\UO, \UD,\MP,\MA$ and uniformly in~$p\in(0,1)$,
\begin{equation}
\label{E:FEapx}
\lim_{L\to\infty}\,F_{L,\alpha}(p)=F_\alpha(p).
\end{equation}
\end{theorem}

\begin{proofsect}{Proof}
This is a result of standard calculations of Gaussian integrals in momentum representation. 
We begin by noting that the Lebesgue measure $\prod_x\textd\phi_x$ can be regarded as the product of~$\nu_L$, acting only on the gradients of~$\phi$, and~$\textd\phi_z$ for some fixed~$z\in\T_L$. Neglecting temporarily the \emph{a priori} bond weights~$p$ and~$(1-p)$, the partition function~$Z_{L,\alpha}$, $\alpha=\Oind,\Dind,\UO, \UD,\MP,\MA$, is thus the integral of the Gaussian weight~$(2\pi|\T_L|)^{-1/2}\texte^{-\frac12(\phi,C_\alpha^{-1}\phi)}$ against the measure~$\prod_x\textd\phi_x$, where the covariance matrix~$C_\alpha$ is defined by the quadratic form
\begin{equation}
\label{phiCphi}
(\phi,C_\alpha^{-1}\phi)=
\sum_{b\in\B_L}\kappa_b^{(\alpha)}(\nabla_b\phi)^2+\frac1{|\T_L|}\Bigl(\sum_{x\in\T_L}\phi_x\Bigr)^2.
\end{equation}
Here~$(\kappa_b^{(\alpha)})$ are the bond weights of pattern~$\alpha$. 
Indeed, the integral over 
$\textd \phi_z$  with the gradient variables fixed yields $(2\pi|\T_L|)^{1/2}$ which cancels the term in front of the Gaussian weight.
The purpose of the above rewrite was to reinsert the ``zero mode'' $\wh\phi_{\bzero}=|\T_L|^{-1/2}\sum_x\phi_x$ into the partition function; $\wh\phi_{\bzero}$ was not subject to integration due to the restriction to gradient variables.

To compute the Gaussian integral, we need to diagonalize~$C_\alpha$. For that we will pass to the Fourier components $\wh\phi_{\bk}=|\T_L|^{-1/2}\sum_{x\in\T_L}\phi_x\texte^{\texti\bk\cdot x}$ with the result
\begin{equation}
\label{E:3.20}
(\phi,C_\alpha^{-1}\phi)=\sum_{\bk,\bk'\in\wt\T_L}
\wh\phi_{\bk}\wh\phi_{\bk'}^*\biggl(\delta_{\bk,\bzero}\delta_{\bk',\bzero}+
\sum_{\sigma=1,2}A_{\bk,\bk'}^{(\sigma)}(1-\texte^{-\texti k_\sigma})(1-\texte^{\texti k'_\sigma})\biggr),
\end{equation}
where $\wt\T_L=\{\frac{2\pi}L(n_1,n_2)\colon 0\le n_1,n_2<L\}$ is the reciprocal torus, ~$\delta_{\bp,\bq}$ is the Kronecker delta and
\begin{equation}
A_{\bk,\bk'}^{(\sigma)}=\frac1{|\T_L|}\sum_{x\in\T_L}\kappa^{(\alpha)}_{(x,x+\hate_\sigma)}\,\texte^{\texti(\bk'-\bk)\cdot x}.
\end{equation}
Now if the horizontal part of~$(\kappa_b^{(\alpha)})$ is translation invariant in the~$\gamma$-th direction, then~$A_{\bk,\bk'}^{(1)}=0$ whenever~$k_\gamma\ne k_\gamma'$, while if it is ``only'' 2-periodic, then~$A_{\bk,\bk'}^{(1)}=0$ unless~$k_\gamma=k_\gamma'$ or~$k_\gamma=k_\gamma'+\pi\mod2\pi$. Similar statements apply to the vertical part of~$(\kappa_b^{(\alpha)})$ and~$A_{\bk,\bk'}^{(2)}$. 
Since all of our partition functions come from 2-periodic configurations, the covariance matrix can be cast into a block-diagonal form, with~$4\times4$ blocks~$\Theta_\alpha(\bk)$ collecting all matrix elements that involve the momenta~$(\bk,\bk+\pi\hate_1,\bk+\pi\hate_2,\bk+\pi\hate_1+\pi\hate_2)$. Due to the reinsertion of the  
``zero mode''---cf~\eqref{phiCphi}---all 
of these blocks are non-singular (see also the explicit calculations below). 

Hence we get that, for all~$\alpha=\Oind,\Dind,\UO,\UD,\MP,\MA$,
\begin{equation}
\label{E:3.22a}
\frac{Z_{L,\alpha}}{(2\pi)^{\frac12(L^2-1)}}=\frac1L\,p^{\NO}(1-p)^{\ND}\prod_{\bk\in\wt\T_L}\biggl[\frac1{\det\Theta_\alpha(\bk)}\biggr]^{\ffrac18},
\end{equation}
where $\NO$ and~$\ND$ denote the numbers of ordered and disordered bonds in the underlying bond configuration and where the exponent~$\ffrac18$ takes care of the fact that in the product, each~$\bk$ gets involved in \emph{four} distinct terms.
Taking logarithms and dividing by~$|\T_L|$, the sum over the reciprocal torus converges to a Riemann integral over the Brillouin zone $[-\pi,\pi]\times[-\pi,\pi]$ (the integrand has only logarithmic singularities in all cases, which are harmless for this limit).

It remains to justify the explicit form of the free energies in all cases under considerations. Here the situations~$\alpha=\Oind,\Dind,\MP$ are fairly standard, so we will focus on~$\alpha=\UO$ and~$\alpha=\MA$ for which some non-trivial calculations are needed. In the former case we get that
\begin{equation}
A_{\bk,\bk'}^{(1)}=\kappaO\delta_{\bk,\bk'},
\quad
A_{\bk,\bk}^{(2)}=\frac{\kappaO+\kappaD}2,
\quad
A_{\bk,\bk+\pi\hate_1}^{(2)}=\frac{\kappaO-\kappaD}2,
\end{equation}
with~$A_{\bk,\bk'}^{(\sigma)}=0$ for all values that are not of this type. Plugging into \eqref{E:3.20} we find that the $(\bk,\bk+\pi\hate_1)$-subblock of~$\Theta_\UO(\bk)$ reduces essentially to the~$2\times2$-matrix in \eqref{PiUO}. Explicitly,
\begin{equation}
\label{E:3.24a}
\Theta_\UO(\bk)=
\text{\rm diag}\bigl(\delta_{\bk,\bzero},\delta_{\bk,\pi\hate_1},\delta_{\bk,\pi\hate_2},\delta_{\bk,\pi\hate_1+\pi\hate_2}\bigr)
+\left(\begin{matrix}
\Pi_\UO(\bk) & 0 \\
0 & \Pi_\UO(\bk+\pi\hate_2)\\
\end{matrix}\right).
\end{equation}
Since~$k'_\sigma=k_\sigma$ whenever~$A_{\bk,\bk'}^{(\sigma)}\ne0$, the block matrix~$\Theta_\UO(\bk)$ will only be a function of moduli-squared of~$a_\pm$ and~$b_-$.
Using \eqref{E:3.24a} in \eqref{E:3.22a} we get~\eqref{E:FUO}. 

As to the~$\MA$-case the only non-zero elements of~$A_{\bk,\bk'}^{(\sigma)}$ are
\begin{equation}
A_{\bk,\bk}^{(1)}=A_{\bk,\bk}^{(2)}=\frac{\kappaO+\kappaD}2
\quad\text{and}\quad
A_{\bk,\bk+\pi\hate_2}^{(1)}=A_{\bk,\bk+\pi\hate_1}^{(2)}=\frac{\kappaO-\kappaD}2.
\end{equation}
So, again, $k'_\sigma=k_\sigma$ whenever~$A_{\bk,\bk'}^{(\sigma)}\ne0$ and so~$\Theta_\MA(\bk)$ depends only on~$|a_\pm|^2$ and~$|b_\pm|^2$. An explicit calculation shows that
\begin{equation}
\Theta_\MA(\bk)=
\text{\rm diag}\bigl(\delta_{\bk,\bzero},\delta_{\bk,\pi\hate_1},\delta_{\bk,\pi\hate_2},\delta_{\bk,\pi\hate_1+\pi\hate_2}\bigr)+
\Bigl(\frac{\kappaO-\kappaD}2\Bigr)\,\Pi_\MA(\bk),
\end{equation}
where~$\Pi_\MA(\bk)$ is as in \eqref{E:Pi}. Plugging into \eqref{E:3.22a}, we get \eqref{E:FMA}.
\end{proofsect}

\subsection{Optimal patterns}
Next we establish the crucial fact that the spin-wave free energies corresponding to inhomogeneous patterns~$\UO,\UD,\MP,\MA$ exceed the smaller of $F_{\Oind}$ and~$F_{\Dind}$ by a quantity that is large, independent of~$p$, once~$\kappaO\gg\kappaD$.

\begin{theorem}
\label{T:min}
There exists  $c_1\in \mathbb R$  such that if $\kappaD\le\xi\,\kappaO$ with $\xi\in(0,1)$,
then for all~$p\in(0,1)$,
\begin{equation}
\label{E:min}
\min_{\alpha=\UO,\UD,\MP,\MA}F_\alpha(p)-\min_{\tilde\alpha=\Oind,\Dind} F_{\tilde\alpha}(p)\ge \frac18\log\frac{\kappaO}{\kappaD}
+\frac14\log(1-\xi)+c_1.
\end{equation}
\end{theorem}

\begin{proof}
Let us use $I$ and~$J$ to denote  the integrals
\begin{equation}
\label{E:I}
I=\frac12\int_{[-\pi,\pi]^2}\frac{\textd\bk}{(2\pi)^2}\log\bigl\{\widehat D(\bk)\bigr\},
\quad\text{and}\quad
J=\int_{[-\pi,\pi]^2}\frac{\textd\bk}{(2\pi)^2}\log\bigl| a_-\bigr|.
\end{equation}
We will prove \eqref{E:min} with $c_1=J-I$.

First, we have
\begin{equation}
F_{\Oind}(p)=-2\log p+\frac12\log\kappaO +I
\end{equation}
and
\begin{equation}
F_{\Dind}(p)=-2\log(1-p)+\frac12\log \kappaD +I,
\end{equation}
while an inspection of \eqref{E:FMA} yields
\begin{multline}
\quad
F_{\MA}(p)=-\log\bigl[p(1-p)\bigr]
\\
+\frac18\int_{[-\pi,\pi]^2}\frac{\textd\bk}{(2\pi)^2}\log\Bigl\{\kappaO\kappaD    (\kappaO-\kappaD)^2
\bigl|a_+ a_- b_+ b_-\bigr|^2\Bigr\}
\\\ge  
-\log\bigl[p(1-p)\bigr]+\frac38{\log \kappaO}+\frac18\log\kappaD
+\frac14\log(1-\xi) +J.
\end{multline}
Using that
\begin{equation}
\min\bigl\{F_{\Oind},F_{\Dind}\bigr\}\le \frac12\bigl( F_{\Oind}+F_{\Dind}\bigr),
\end{equation}
we thus get
\begin{equation}
F_{\MA}(p)-\min\bigl\{F_{\Oind}(p),F_{\Dind}(p)\bigr\}\ge
\frac18 \log\frac{\kappaO}{\kappaD} +\frac14\log (1-\xi) +J-I,
\end{equation}
which agrees with \eqref{E:min} for our choice of~$c_1$.

Coming to the free energy $F_{\UO}$, using \eqref{detUO-bound} we evaluate
\begin{equation}
\det \Pi_{\UO}(\bk)\ge  \kappaO^2 |a_-|^2 |a_+|^2
=\Bigl(\frac{\kappaO}{\kappaD}\Bigr)^{\ffrac12}  \kappaO^{\ffrac32} \kappaD^{1/2}   |a_-|^2 |a_+|^2
\end{equation}
yielding
\begin{equation}
F_{\UO}(p)\ge -\frac12\log\bigl[p^3(1-p)\bigr]+\frac18\log\frac{\kappaO}{\kappaD}
+\frac38{\log \kappaO}+\frac18{\log \kappaD}+J.
\end{equation}
Bounding
\begin{equation}
\min\bigl\{F_{\Oind},F_{\Dind}\bigr\}\le \frac34 F_{\Oind}+\frac14 F_{\Dind}
\end{equation}
we thus get
\begin{equation}
\label{E:FUObd}
F_{\UO}(p)-\min\bigl\{F_{\Oind}(p),F_{\Dind}(p)\bigr\}\ge
\frac18  \log\frac{\kappaO}{\kappaD} +  J-  I,
\end{equation}
in agreement with \eqref{E:min}.
The computation for $F_{\UD}$ is completely analogous, interchanging only the roles of~$\kappaO$ and~$\kappaD$ as well as~$p$ and~$(1-p)$.
From the lower bound
\begin{equation}
\det \Pi_{\UO}(\bk) \ge  \kappaO \kappaD |b_-|^4
=\Bigl(\frac{\kappaO}{\kappaD}\Bigr)^{1/2}  \kappaO^{1/2} \kappaD^{3/2}   |b_-|^4
\end{equation}
and the inequality
\begin{equation}
\min\bigl\{F_{\Oind},F_{\Dind}\bigr\}\le \frac14 F_{\Oind}+\frac34 F_{\Dind},
\end{equation}
we get again
\begin{equation}
F_{\UD}(p)-\min\bigl\{F_{\Oind}(p),F_{\Dind}(p)\bigr\}\ge
\frac18  \log\frac{\kappaO}{\kappaD} +  J-  I,
\end{equation}
which is identical to \eqref{E:FUObd}.

Finally, for the free energy $F_{\MP}$, we first note that
\begin{equation}
F_{\MP}(p)
\ge -\log\bigl[p(1-p)\bigr]
+\frac12\log\kappaO +J,
\end{equation}
which yields
\begin{equation}
F_{\MP}(p)-F_{\Dind}(p) \ge \log\frac{1-p}p
+\frac12\log\frac{\kappaO}{\kappaD} +J-I.
\end{equation}
Under the condition that~$\log\frac{1-p}p\ge-\frac38\log\frac{\kappaO}{\kappaD}$, we again get~\eqref{E:min}.
For the complementary values of~$p$, we will compare~$F_{\MP}$ with~$F_{\Oind}$:
\begin{equation}
F_{\MP}(p)-F_{\Oind}(p) \ge \log\frac p{1-p}+J-I,
\end{equation}
Since we now have~$\log\frac{1-p}p\le-\frac38\log\frac{\kappaO}{\kappaD}$, this yields \eqref{E:min} with the above choice of~$c_1$.
\end{proof}

\section{Proof of phase coexistence}
\label{sec:PT}\noindent
In this section we will apply the calculations from the previous section to the proof of Theorems~\ref{T:main} and~\ref{T:torus}. Throughout this section we assume that~$\kappaO>\kappaD$
and that~$L$ is even. 
We begin with a review of the technique of chessboard estimates
which, for later convenience, we formulate directly in terms of extended configurations~$(\kappa_b,\eta_b)$.

\subsection{Review of RP/CE technology}
\label{S:chessboard}\noindent
Our principal tool will be chessboard estimates, based on reflection positivity. 
To define these concepts, let us consider the torus $\T_L$
and let us split~$\T_L$ into two symmetric halves, $\T_L^+$ and $\T_L^-$, sharing a 
``plane of sites''
on their boundary.
We will refer to the set $\T_L^+\cap\T_L^-$ as \emph{plane of reflection} and denote it by $P$.
The half-tori $\T_L^\pm$ inherit the nearest-neighbor structure from $\T_L$; 
we will use $\B_L^\pm$ to denote the corresponding sets of edges.
On the extended configuration space, there is a canonical map $\theta_P\colon\R^{\B_L}\times\{\kappaO,\kappaD\}^{\B_L}\to\R^{\B_L}\times\{\kappaO,\kappaD\}^{\B_L}$---induced by the reflection
of $\T_L^+$ into $\T_L^-$ through~$P$---which is defined as follows:
If $b,b'\in\B_L$ are related via~$b'=\theta_P(b)$, then we put
\begin{equation}
\label{E:vartheta}
\bigl(\theta_P\eta\bigr)_b=
\begin{cases}
-\eta_{b'},\qquad&\text{if }b\perp P,
\\
\eta_{b'},\qquad&\text{if }b\parallel P,
\end{cases}
\end{equation}
and
\begin{equation}
\label{E:vartheta2}
(\theta_P\kappa)_b=\kappa_{b'}.
\end{equation}
Here~$b\perp P$ denotes that~$b$ is orthogonal to~$p$ while~$b\parallel P$ indicates that~$b$ is parallel to~$P$.
The minus sign in the case when~$b\perp P$ is fairly natural if we recall that~$\eta_b$ represents the difference of~$\phi_x$ between the endpoints of~$P$. This difference changes sign under reflection through~$P$ if~$b\perp P$ and does not if~$b\parallel P$.

\smallskip
Let $\scrF_P^\pm$ be the $\sigma$-algebras of events that depend only on the portion of $(\eta_b,\kappa_b)$-configuration on~$\B_L^\pm$;
explicitly  $\scrF_P^\pm=\sigma\bigl(\eta_b,\kappa_b; b\in\B_L^\pm\bigr)$.
Reflection positivity is, in its essence, a bound on the correlation between events (and random variables) from $\scrF_P^+$ and $\scrF_P^-$.
The precise definition is as follows:

\begin{mydefinition}
\label{D:RP}
Let $\BbbP$ be a probability measure on configurations $(\eta_b,\kappa_b)_{b\in\B_L}$ and let $\E$ be the corresponding expectation.
We say that~$\BbbP$ is \emph{reflection positive} if for any plane of reflection~$P$ and any two bounded $\scrF_P^+$-measurable random variables
$X$ and $Y$ the following inequalities hold:
\begin{equation}
\label{E:XY}
\E(X\theta_P(Y))=\E(Y\theta_P(X))
\end{equation}
and
\begin{equation}
\label{E:X}
\E(X\theta_P(X))\ge 0.
\end{equation}
Here, $\theta_P(X)$ denotes the random variable $X\circ \theta_P$.
\end{mydefinition}

Next we will discuss how reflection positivity underlines our principal technical tool: chessboard estimates.
Consider an event $\AA$ that depends only on the $(\eta_b,\kappa_b)$-configurations on the plaquette with the lower-left corner at the torus origin.
We will call such an~$\AA$ a \emph{plaquette event}.
For each $x\in\T_L$, we define $\vartheta_x(\AA)$ to be the event depending only on the configuration on the plaquette 
with the lower-left corner at $x$ which is obtained from~$\AA$ as follows:
If both components of $x$ are even, then  $\vartheta_x(\AA)$ is simply the translate of $\AA$ by $x$.
In the remaining cases we first reflect $\AA$ along the side(s) of the plaquette in the direction(s) where the component of $x$ is odd,
and then translate the resulting event appropriately. (Thus, there are four possible ``versions'' of $\vartheta_x(\AA)$,
depending on the parity of $x$.)

\smallskip
Here is the desired consequence of reflection positivity:

\begin{theorem}[Chessboard estimate]
\label{T:chess}
Let $\BbbP$ be a reflection-positive measure on configurations $(\eta_b,\kappa_b)_{b\in\B_L}$. Then for any plaquette events $\AA_1,\dots,\AA_m$ and any distinct sites $x_1,\dots,x_m\in\T_L$, 
\begin{equation}
\label{E:chess}
\BbbP\Bigl(\,\bigcap_{j=1}^m\vartheta_{x_j}(\AA_j)\Bigr)\le
\prod_{j=1}^m\,\BbbP\Bigl(\,\bigcap_{x\in\T_L}\vartheta_{x}(\AA_j)\Bigr)^{\frac1{|\T_L|}}.
\end{equation}
\end{theorem}
\begin{proof}
See \cite[Theorem~2.2]{FL}.
\end{proof}

The moral of this result---whose proof boils down to the Cauchy-Schwarz inequality for the inner product
$X,Y\mapsto \E(X\theta_P(Y))$---is that the probability of any number of plaquette events factorizes, as a bound, into the product of probabilities. 
This is particularly useful for contour estimates (of course, provided that the word contour refers to a collection of plaquettes on each of which some ``bad'' event occurs). 
Indeed, by \eqref{E:chess} the probability of a contour will be suppressed exponentially in the number of constituting plaquettes.
 
In light of \eqref{E:chess}, our estimates will require good bounds on probabilities of the so called \emph{disseminated events} $\bigcap_{x\in\T_L}\vartheta_{x}(\AA)$. 
Unfortunately, the event $\AA$ is often a conglomerate of several,
more elementary events which makes a direct estimate of $\bigcap_{x\in\T_L}\vartheta_{x}(\AA)$ complicated.
Here the following subadditivity property will turn out to be useful.

\begin{mylemma}[Subadditivity]
\label{L:sub}
Suppose that~$\BbbP$ is a reflection-positive measure and let $\AA_1,\AA_2,\dots$ and $\AA$ be plaquette events such that $\AA\subset \bigcup_j\AA_j$. Then
\begin{equation}
\label{E:sub}
\BbbP\Bigl(\bigcap_{x\in\T_L}\vartheta_{x}(\AA)\Bigr)^{\frac1{|\T_L|}}\le
\sum_j\,\BbbP\Bigl(\bigcap_{x\in\T_L}\vartheta_{x}(\AA_j)\Bigr)^{\frac1{|\T_L|}}.
\end{equation}
\end{mylemma}
\begin{proof}
This is Lemma 6.3 of~\cite{BCN}.
\end{proof}

Apart from the above reflections, which we will call \emph{direct}, one estimate---namely \eqref{diag-refl}---in the proof of Theorem~\ref{T:main} requires the use of so called \emph{diagonal reflections}. Assuming~$L$ is even, these are reflections in the planes~$P$ of sites of the form
\begin{equation}
P_\pm=\bigl\{y\in\T_L\colon\hate_1\cdot(y-x)\mp\hate_2\cdot(y-x)\in\{0,\ffrac L2\}\bigr\}.
\end{equation}
Here~$x$ is a site that the plane passes through and~$\hate_1$ and~$\hate_2$ are the unit vectors in the $x$ and~$y$-coordinate directions. As before, the plane has two components---one corresponding to $\hate_1\cdot(y-x)=\pm\hate_2\cdot(y-x)$ and the other corresponding to $\hate_1\cdot(y-x)=\pm\hate_2\cdot(y-x)+\ffrac L2$---and it divides~$\T_L$ into two equal parts. This puts us into the setting assumed in Definition~\ref{D:RP}. Some care is needed in the definition of reflected configurations: If~$b'$ is the bond obtained by reflecting~$b$ through~$P$, then
\begin{equation}
\label{E:vartheta-diag}
\bigl(\theta_P\eta\bigr)_b=
\begin{cases}
\eta_{b'},\qquad&\text{if }P=P_+,
\\
-\eta_{b'},\qquad&\text{if }P=P_-.
\end{cases}
\end{equation}
This is different compared to \eqref{E:vartheta} because the reflection in~$P_+$ preserves orientations of the edges, while that in~$P_-$ reverses them. 

\begin{myremark}
While we will only apply these reflections in~$d=2$, we note that the generalization to higher dimensions is straightforward; just consider all planes as above with~$(\hate_1,\hate_2)$ replaced by various pairs~$(\hate_i,\hate_j)$ of distinct coordinate vectors. These reflections will of course preserve the orientations of all edges in directions distinct from~$\hate_i$ and~$\hate_j$.
\end{myremark}

\subsection{Phase transitions on tori}
\label{sec4.2}\noindent
Here we will provide the proof of phase transition in the form stated in Theorem~\ref{T:torus}. We follow pretty much the standard approach to proofs of order-disorder transitions which dates all the way back to~\cite{DS,KS,KS-proceedings}. A somewhat different approach (motivated by another perspective) to this proof can be found in~\cite{BK05}.

\smallskip
In order to use the techniques decribed in the previous section, we have to 
determine when the extended gradient Gibbs measure $Q_L$ on~$\T_L$ obeys the conditions of reflection positivity.

\begin{myproposition}
\label{P:RP}
Let $V$ be of the form \eqref{V-rep}
with any probability measure~$\varrho$ for which~$Z_{L,V}<\infty$.
Then~$Q_L$ is reflection positive for both direct and diagonal reflections.
\end{myproposition}

\begin{proof}
The proof is the same for both types of reflections so we we proceed fairly generally. Pick a plane of reflection~$P$. Let~$z$ be a site on~$P$ and let us reexpress the~$\eta_b$'s back in terms of the~$\phi$'s with the convention that~$\phi_z=0$. Then
\begin{equation}
\nu_L(\textd\eta_b)=\delta(\textd\phi_z)\prod_{x\ne z}\textd\phi_x.
\end{equation}
Next, let us introduce the quantity
\begin{equation}
W(\eta,\kappa)=\frac12\sum_{b\in\B_L^+\smallsetminus P}\kappa_b\eta_b^2+\frac14\sum_{b\in P}\kappa_b\eta_b^2.
\end{equation}
(We note in passing that the removal of~$P$ from the first sum is non-trivial even for diagonal reflections once $d\ge3$.)
Clearly,~$W$ is~$\scrF_P^+$-measurable and the full $(\eta,\kappa)$-interaction is simply~$W(\eta,\kappa)+(\theta_PW)(\eta,\kappa)$. The Gibbs measure~$Q_L$ can then be written
\begin{equation}
Q_L(\textd\eta,\textd\kappa)=\frac1{Z_{L,V}}\texte^{-W(\nabla\phi,\kappa)-(\theta_PW)(\nabla\phi,\kappa)}
\delta(\textd\phi_z)\biggl(\,\prod_{x\ne z}\textd\phi_x\biggr)\prod_{b\in\B_L}\rho(\textd\kappa_b)
\end{equation}
Now pick a bounded, $\scrF_P^+$-measurable function~$X=X(\eta,\kappa)$ and integrate the function~$X\,\theta_P X$ with respect to the torus measure~$Q_L$. 
If~$\scrG_P$ is the $\sigma$-algebra generated by random variables~$\phi_x$ and~$\kappa_b$, with~$x$ and~$b$ ``on''~$P$, we have
\begin{equation}
E_{Q_L}\bigl(X\theta_PX\big|\scrG_P)\propto
\biggl(\,\int X(\nabla\phi,\kappa)\texte^{-W(\nabla\phi,\kappa)}\!\!
\prod_{x\in\T_L^+\smallsetminus P}\!\!\textd\phi_x
\prod_{b\in\B_L^+\smallsetminus P}\!\!\varrho(\textd\kappa_b)\biggr)^2\ge0,
\end{equation}
where the values of~$(\kappa_b,\phi_x)$ on~$P$ are implicit in the integral.
This proves the property in~\eqref{E:X}; the identity \eqref{E:XY} follows by the reflection symmetry of~$Q_L$.
\end{proof}

Let us consider two good plaquette events,~$\GG_\ord$ and~$\GG_\dis$, that all edges on the plaquette are ordered and disordered, respectively. Let~$\BB=(\GG_\ord\cup\GG_\dis)^\cc$ denote the corresponding bad event. Given a plaquette event~$\AA$, let
\begin{equation}
\zz_{L,p}(\AA)=\biggl[Q_L\Bigl(\,\bigcap_{x\in\T_L}\vartheta_{x}(\AA)\Bigr)\biggr]^{\frac1{|\T_L|}}
\end{equation}
abbreviate the quantity on the right-hand side of \eqref{E:chess} and define
\begin{equation}
\label{cistezz}
\zz(\AA)=\limsup_{L\to\infty}\sup_{0\le p\le 1}\zz_{L,p}(\AA).
\end{equation}
The calculations from Sect.~\ref{S:spin-wave} then permit us to draw the following conclusion:

\begin{mylemma}
\label{L:bad}
For each~$\delta>0$ there exists~$c>0$ such that if~$\kappaO\ge c\kappaD$, then
\begin{equation}
\label{BB-bd}
\zz(\BB)<\delta.
\end{equation}
Moreover, there exist~$p_0,p_1\in(0,1)$ such that
\begin{equation}
\label{dis-bd}
\limsup_{L\to\infty}\zz_{L,p}(\GG_\ord)<\delta,\qquad p<p_0,
\end{equation}
and
\begin{equation}
\label{ord-bd}
\limsup_{L\to\infty}\zz_{L,p}(
\GG_\dis)<\delta,\qquad p>p_1.
\end{equation}
\end{mylemma}

\begin{proofsect}{Proof}
The event~$\BB$ can be decomposed into a disjoint union of events~$\BB_i$ each of which admits exactly one arrangement of ordered and disordered bonds around the plaquette; see Fig.~\ref{fig3} for the relevant patterns. If~$\BB_i$ is an event of type~$\alpha\in\{\Oind,\Dind,\UO, \UD,\MP,\MA\}$, then
\begin{equation}
\limsup_{L\to\infty}\zz_{L,p}(\BB_i)
\le\exp\bigl\{-[F_{L,\alpha}(p)-\min_{\tilde\alpha=\Oind,\Dind} F_{\tilde\alpha}(p)]\bigr\}.
\end{equation}
By Theorem~\ref{T:min}, the right-hand side is bounded
by $O(1)(\ffrac\kappaO\kappaD)^{-1/8}$, uniformly in~$p$.
Applying Lemma~\ref{L:sub}, we conclude that~$\zz_{L,p}(\BB)$ is small uniformly in~$p\in[0,1]$ once~$L\gg1$.
(The values~$p=0,1$ are handled by a limiting argument.)

The bounds \twoeqref{dis-bd}{ord-bd} follow by the fact that
\begin{equation}
F_{\Oind}(p)-F_{\Dind}(p)=-2\log\frac p{1-p}+\log\frac\kappaO\kappaD,
\end{equation}
which is (large) negative for~$p$ near one and (large) positive for~$p$ near zero.
\end{proofsect}

From~$\zz(\BB)\ll1$ we immediately infer that the bad events occur with very low frequency. Moreover, a standard argument shows that the two good events do not like to occur in the same configuration. An explicit form of this statement is as follows:

\begin{mylemma}
\label{L:good-good}
Let~$R_L^\ord$ be the random variable from \eqref{Rord}. There 
exists a constant
$C<\infty$ such that for all (even)~$L\ge1$ and all~$p\in[0,1]$,
\begin{equation}
\label{R-cov}
E_{Q_L}\bigl(R_L^\ord(1-R_L^\ord)\bigr)
\le C\zz_{L,p}(\BB).
\end{equation}
\end{mylemma}

\begin{proofsect}{Proof}
The claim follows from the fact that, for some constant~$C'<\infty$,
\begin{equation}
\label{E:ord-dis}
Q_L\bigl(\vartheta_x(\GG_\ord)\cap\vartheta_y(\GG_\dis)\bigr)\le C'\zz_{L,p}(\BB)^4,
\end{equation}
uniformly in~$x,y\in\T_L$.
Indeed, the expectation in \eqref{R-cov} is the average of the probabilities~$Q_L(\kappa_b=\kappaO,\,\kappa_{\tilde b}=\kappaD)$ over all~$b,\tilde b\in\B_L$. If~$x$ and~$y$ denotes the plaquettes containing the bonds~$b$ and~$\tilde b$, respectively, then this probability is bounded by~$Q_L(\vartheta_x(\GG_\ord)\cap\vartheta_y(\GG_\ord^\cc))$. But~$\GG_\ord^\cc=\BB\cup\GG_\dis$ and so by \eqref{E:ord-dis} the latter probability is bounded by $\zz_{L,p}(\BB)+C'\zz_{L,p}(\BB)^4\le(C'+1)\zz_{L,p}(\BB)$, where we used~$\zz_{L,p}(\BB)\le1$.

It remains to prove \eqref{E:ord-dis}. Consider the event~$\vartheta_x(\GG_\ord)\cap\vartheta_y(\GG_\dis)$ where, without loss of generality,~$x\ne y$. We claim that on this event, the good plaquettes at~$x$ and~$y$ are separated from each other by a $*$-connected circuit of bad plaquettes. To see this, consider the largest connected component of good plaquettes containing~$x$ and note that no plaquette neighboring on this component can be good, because (by definition) the events~$\GG_\ord$ and~$\GG_\dis$ cannot occur at neighboring plaquettes (we are assuming that~$\kappaO\ne\kappaD$). By chessboard estimates, the probability in~$Q_L$ of any such (given) circuit is bounded by~$\zz_{L,p}(\BB)$ to its size; a standard Peierls' argument in toroidal geometry (cf the proof of~\cite[Lemma~3.2]{BCN}) now shows that the probability in \eqref{E:ord-dis} is dominated by the probability of the shortest possible contour---which is~$\zz_{L,p}(\BB)^4$.
(The contour argument requires that~$\zz_{L,p}(\BB)$ be smaller than some constant, but this we may assume to be automatically satisfied because the left-hand side of \eqref{R-cov} is less than one.)
\end{proofsect}

Now we are in a position to prove our claims concerning the torus state:

\begin{proofsect}{Proof of Theorem~\ref{T:torus}}
Let~$R^\ord_L$ be the fraction of ordered bonds on~$\T_L$ (cf. \eqref{Rord}) and let~$\chi_L(p)$ be the expectation of~$R_L^\ord$ in the extended torus state~$Q_L$ with parameter~$p$. Since $(1-p)^{-|\B_L|}Z_{L,V}$ is log-convex in the variable~$h=\log\frac p{1-p}$, 
and
\begin{equation}
|\B_L | \chi_L(p) =\frac{\partial \log\bigl((1-p)^{-|\B_L|}Z_{L,V}\bigr)}{\partial h} ,
\end{equation}
we can conclude that
the function $p\mapsto\chi_L(p)$ is non-decreasing. Moreover, as the thermodynamic limit of the torus free energy exists (cf Proposition~\ref{prop-FE} in Sect.~\ref{sec5.3}), the limit $\chi(p)=\lim_{L\to\infty}\chi_L(p)$ exists at all but perhaps a countable number of~$p$'s---namely the set
$\DD\subset[0,1]$ of points where the limiting free energy is not differentiable.

Next we claim that $Q_L(|R_L^\ord-\chi(p)|>\epsilon_0)$ tends to zero as~$L\to\infty$ for
all~$\epsilon_0>0$ and all~$p\not\in\DD$. Indeed, if this probability 
stays
uniformly positive along some subsequence of~$L$'s for 
some~$\epsilon_0>0$,
then the boundedness of
$R_L^\ord$ ensures that for some~$\zeta>0$ and some~$\epsilon>0$ we have
$Q_L(R_L^\ord>\chi(p)+\epsilon)\ge\zeta$ \emph{and}
$Q_L(R_L^\ord<\chi(p)-\epsilon)\ge\zeta$ for all~$L$ in this 
subsequence.
Vaguely speaking, this implies~$p\in\DD$ because one is then able to extract two infinite-volume Gibbs states with distinct densities of ordered bonds. A formal proof goes as follows: Consider the cumulant generating function~$\Psi_L(h)=|\B_L|^{-1}\log E_{Q_L}(\texte^{h|B_L|R_L^\ord})$ and note that its thermodynamic limit, $\Psi(h)=\lim_{L\to\infty}\Psi_L(h)$, is convex in~$h$ and differentiable at~$h=0$ whenever~$p\not\in\DD$. But~$Q_L(R_L^\ord>\chi(p)+\epsilon)\ge\zeta$ in conjunction with the exponential Chebyshev inequality implies
\begin{equation}
\Psi_L(h)-h\bigl(\chi(p)+\epsilon\bigr)\ge\frac{\log\zeta}{|\B_L|},
\end{equation}
which by taking~$L\to\infty$ and~$h\downarrow0$ yields a lower bound on the right derivative at origin,  $\frac{\textd}{\textd h_+}\Psi(h)\ge\chi(p)+\epsilon$. By the same token $Q_L(R_L^\ord<\chi(p)-\epsilon)\ge\zeta$ implies 
an upper bound on the left derivative, $\frac{\textd}{\textd h_-}\Psi(h)\le\chi(p)-\epsilon$. Hence, both probabilities can be uniformly positive only if~$p\in\DD$.

To prove the desired claim it remains to show that~$\chi$ jumps from values near zero to values near one at some~$\pt\in(0,1)$.
To this end we first observe that
\begin{equation}
\lim_{L\to\infty}E_{Q_L}\bigl(R_L^\ord(1-R_L^\ord)\bigr)=\chi(p)\bigl[1-\chi(p)\bigr],
\qquad p\not\in\DD.
\end{equation}
This follows by the fact that on the event~$\{\chi(p)-\epsilon<R_L^\ord<\chi(p)+\epsilon\}$---whose probability tends to one as~$L\to\infty$---the quantity~$R_L^\ord(1-R_L^\ord)$ is bounded between~$[\chi(p)+\epsilon](1-\chi(p)+\epsilon)$ and $[\chi(p)-\epsilon](1-\chi(p)-\epsilon)$
provided~$\epsilon\le\min\{\chi(p),1-\chi(p)\}$. 
Lemma~\ref{L:good-good} now implies
\begin{equation}
\chi(p)\bigl[1-\chi(p)\bigr]\le C\zz(\BB),
\end{equation}
with~$\zz(\BB)$ defined in \eqref{cistezz}. 
By Lemma~\ref{L:bad}, for each~$\delta>0$ there is a constant~$c>0$ such that 
\begin{equation}
\chi(p)\in[0,\delta]\cup[1-\delta,1],\qquad p\not\in\DD,
\end{equation}
once~$\ffrac\kappaO\kappaD\ge c$.
But the bounds \twoeqref{dis-bd}{ord-bd} ensure that~$\chi(p)\in[0,\delta]$ for $p\ll1$ and~$\chi(p)\in[1-\delta,1]$ for $1-p\ll1$. Hence, by the monotonicity of~$p\mapsto\chi(p)$, there exists a unique value~$\pt\in(0,1)$ such that~$\chi(p)\le\delta$ for~$p<\pt$ while $\chi(p)\ge1-\delta$ for~$p>\pt$. In light of our previous reasoning, this proves the bounds \twoeqref{below-pt}{above-pt}.
\end{proofsect}

\subsection{Phase coexistence in infinite volume}
\label{sec4.3}\noindent
In order to prove Theorem~\ref{T:main}, we will need to derive a concentration bound on the tilt of the torus states. This is the content of the following lemma:

\begin{mylemma}
\label{L:tilt}
Let~$\Lambda\subset\T_L$ and let~$\B_\Lambda$ be the set of bonds with both ends in~$\Lambda$.
Given a configuration~$(\eta_b)_{b\in\B_L}$, we use~$U_\Lambda=U_\Lambda(\eta)$ to denote the vector
\begin{equation}
U_\Lambda=\Bigl(\,\frac1{|\B_\Lambda|}\sum_{b\in\Bhor_L\cap\B_\Lambda}\!\eta_b\,,
\frac1{|\B_\Lambda|}\sum_{b\in\Bvert_L\cap\B_\Lambda}\!\eta_b\Bigr)
\end{equation}
of empirical tilt of the configuration~$\eta_b$ in~$\Lambda$. 
Suppose that $\kappa_\mini=\inf\supp\varrho>0$.
Then 
\begin{equation}
\label{E:tight-tilt}
P_L\bigl(|U_\Lambda|\ge\delta\bigr)\le4\texte^{-\frac18\kappa_\mini\delta^2 |\B_\Lambda|}
\end{equation}
for each~$\delta>0$, each~$\Lambda\subset\T_L$ and each~$L$.
\end{mylemma}

\begin{proof}
We will derive a bound on the exponential moment of~$U_\Lambda$. 
Let us fix a collection of numbers $(v_b)_{b\in\B_L}\in\R^{|\B_L|}$ and let~$Q_{L,(\kappa_b)}$ be the conditional law of the~$\eta$'s given a configuration of the~$\kappa$'s.
Let
$Q_{L,\mini}$ be the corresponding law when all~$\kappa_b=\kappa_\mini$. 
In view of the fact that~$Q_{L,(\kappa_b)}$ and $Q_{L,\mini}$ are Gaussian measures and~$\kappa_b\ge\kappa_\mini$, we have
\begin{equation}
\Var_{Q_{L,(\kappa_b)}}\biggl(\,\sum_{b\in\B_L}v_b\eta_b\biggr)
\le
\Var_{Q_{L,\mini}}\biggl(\,\sum_{b\in\B_L}v_b\eta_b\biggr).
\end{equation}
(Note that both measures enforce the same loop conditions.)
The right-hand side is best calculated in terms of the gradients. The result is
\begin{equation}
\Var_{Q_{L,(\kappa_b)}}\biggl(\,\sum_{b\in\B_L}v_b\eta_b\biggr)
\le\frac1{\kappa_\mini}\sum_{b\in\B_L}v_b^2.
\end{equation}
The fact that $E_{Q_{L,(\kappa_b)}}(\eta_b)=0$ and the identity
$E(\texte^X)=\texte^{EX+\frac12\Var(X)}$, valid for any Gaussian random variable, now allow us to conclude
\begin{equation}
\label{E:4.22}
E_{Q_L}\biggl(\,\exp\Bigl\{\sum_{b\in\B_L}v_b\eta_b\Bigr\}\biggr)
\le\exp\Bigl\{\frac1{2\kappa_\mini}\sum_{b\in\B_L}v_b^2\Bigr\}.
\end{equation}
Choosing~$v_b=\lambda\cdot b/|\B_\Lambda|$ on~$\B_\Lambda$ and zero otherwise, we get
\begin{equation}
E_{Q_L}(\texte^{\lambda\cdot U_\Lambda})\le\exp\Bigl\{\frac1{2\kappa_\mini}\frac{|\lambda|^2}{|\B_\Lambda|}\Bigr\}.
\end{equation}
Noting that~$|U_\Lambda|>\delta$ implies that at least one of the components of~$U_\Lambda$ is larger (in absolute value) than~$\ffrac\delta2$, the desired bound follows by a standard exponetial-Chebyshev estimate.
\end{proof}

\begin{myremark}
\label{rem4.9}
We note that the symmetry of the law of the~$\eta$'s in~$Q_{L,(\kappa_b)}$ is crucial for the above argument. In particular, it is not clear how to control the tightness of the empirical tilt $U_\Lambda$ in the measure obtained by normalizing $\exp\{-\sum_b V(\eta_b+h_b)\}\nu_L(\textd \eta)$, where~$h_b=h\cdot b$ is a ``built-in'' tilt. In the strictly convex cases, these measures were used by Funaki and Spohn~\cite{Funaki-Spohn} to construct an infinite-volume gradient Gibbs state with a given value of the tilt.
\end{myremark}

\begin{proofsect}{Proof of Theorem~\ref{T:main}}
With Theorem~\ref{T:torus} at our disposal, the argument is fairly straightforward. Consider a weak (subsequential) limit of the torus states at~$p>\pt$ and then consider another weak limit of these states as~$p\downarrow\pt$.
Denote the result by~$\tilde\mu_\ord$. Next let us perform a similar limit as~$p\uparrow\pt$ and let us denote the resulting measure by~$\tilde\mu_\dis$. As is easy to check, both measures are extended gradient Gibbs measures at parameter~$\pt$.

Next we will show that the two measures are distinct measures of zero tilt. To this end we recall that, by \eqref{above-pt} and the invariance of~$Q_L$ under rotations, $\liminf_{L\to\infty}Q_L(\kappa_b=\kappaO)\ge1-\epsilon$ when~$p>\pt$ while \eqref{below-pt} implies that $\limsup_{L\to\infty}Q_L(\kappa_b=\kappaO)\le\epsilon$ when~$p<\pt$. But $\{\kappa_b=\kappaO\}$ is a local event and so
\begin{equation}
\tilde\mu_\ord(\kappa_b=\kappaO)\ge1-\epsilon
\end{equation}
while
\begin{equation}
\tilde\mu_\dis(\kappa_b=\kappaO)\le\epsilon,
\end{equation}
for all~$b$; i.e.,~$\tilde\mu_\ord\ne\tilde\mu_\dis$. Moreover, the bound \eqref{E:tight-tilt}---being uniform in~$p$ and~$L$---survives the above limits unscathed and so the tilt is exponentially tight in volume for both~$\tilde\mu_\ord$ and~$\tilde\mu_\dis$. It follows that~$U_\Lambda\to0$ as~$\Lambda\uparrow\Z^2$ almost surely with respect to both~$\tilde\mu_\ord$ and~$\tilde\mu_\dis$; i.e., both measures are supported entirely on configurations with zero tilt.

It remains to prove the inequalities \twoeqref{E:ener-bd1}{E:ener-bd2} and thereby ensure that the~$\eta$-marginals ~$\mu_\ord$ and~$\mu_\dis$ of~$\tilde\mu_\ord$ and~$\tilde\mu_\dis$, respectively, are distinct as claimed in the statement of the theorem. The first bound is a consequence of the identity
\begin{equation}
\label{E:var-bd}
\lim_{L\to\infty} E_{Q_L}(\kappa_b\eta_b^2)=\frac14,
\end{equation}
which extends via the aforementioned limits to~$\tilde\mu_\ord$ (as well as~$\tilde\mu_\dis$). Indeed, using Chebyshev's inequality and the fact that~$\mu_\ord(\kappa_b=\kappaD)\le\epsilon$ we get
\begin{equation}
\begin{aligned}
\tilde\mu_\ord(\kappaO\eta_b^2\ge\lambda^2)-\epsilon
&\le\tilde\mu_\ord(\kappaO\eta_b^2\ge\lambda^2,\,\kappa_b=\kappaO)
\\
&\le\tilde\mu_\ord(\kappa_b\eta_b^2\ge\lambda^2)\le\frac{E_{\tilde\mu_\ord}(\kappa_b\eta_b^2)}{\lambda^2}=\frac1{4\lambda^2}.
\end{aligned}
\end{equation}
To prove~\eqref{E:var-bd}, the translation 
and rotation
invariance of~$Q_L$ gives us
\begin{equation}
E_{Q_L}(\kappa_b\eta_b^2)=E_{Q_L}\biggl(E_{Q_{L,(\kappa_b)}}\Bigl(\,\frac1{|\B_L|}\sum_{b\in\B_L}\kappa_b\eta_b^2\Bigr)\biggr).
\end{equation}
Let~$Z_{L,(\kappa_b)}$ denote the integral of~$\exp\{-\frac12\sum_b\kappa_b\eta_b^2\}$ with respect to~$\nu_L$.
Since we have $\nu_L(\beta\textd\eta)=\beta^{|\T_L|-1}\nu_L(\textd\eta)$, simple scaling of all fields yields $Z_{L,(\beta\kappa_b)}=
\beta^{-\frac12(|\T_L|-1)}
Z_{L,(\kappa_b)}$. Intepreting the inner expectation above as the (negative) $\beta$-derivative of~$|\B_L|^{-1}\log Z_{L,(\beta\kappa_b)}$ at~$\beta=1$, we get
\begin{equation}
E_{Q_{L,(\kappa_b)}}\Bigl(\,\frac1{|\B_L|}\sum_{b\in\B_L}\kappa_b\eta_b^2\Bigr)=
\frac{|\T_L|-1}{2|\B_L|}.
\end{equation}
From here~\eqref{E:var-bd} follows by taking~$L\to\infty$ on the right-hand side.

As to the inequality \eqref{E:ener-bd2} for the disordered state, here we first use that the diagonal reflection allows us to disseminate the event~$\{\kappa_b\eta_b^2\le\lambda^2\}$ around any plaquette containing~$b$. Explicitly, if~$(b_1,b_2,b_3,b_4)$ is a plaquette, then
\begin{equation}
\label{diag-refl}
Q_L(\kappa_{b_1}\eta_{b_1}^2\le\lambda^2)\le Q_L\Bigl(\,\bigcap_{b=b_1,\dots,b_4}\{\kappa_{b}\eta_{b}^2\le\lambda^2\}\Bigr)^{\ffrac14}.
\end{equation}
(We are using that the event in question is even in~$\eta$ and so the changes of sign of~$\eta_b$ are immaterial.) Direct reflections now permit us to disseminate 
the resulting plaquette event
all over the torus:
\begin{equation}
Q_L(\kappa_b\eta_b^2\le\lambda^2)\le Q_L\Bigl(\,\bigcap_{\tilde b\in\B_L}\{\kappa_{\tilde b}\eta_{\tilde b}^2\le\lambda^2\}\Bigr)^{\frac1{4|\T_L|}}.
\end{equation}
Bounding the indicator of the giant intersection by
\begin{equation}
\texte^{(\beta-1)\lambda^2|B_L|}\exp\Bigl\{-\frac12(\beta-1)\sum_{b\in\B_L}\kappa_b\eta_b^2\Bigr\},
\end{equation}
for~$\beta\ge1$, and invoking the scaling of the partition function~$Z_{L,(\beta\kappa_b)}$, we deduce
\begin{equation}
Q_L(\kappa_b\eta_b^2\le\lambda^2)\le\biggl[\frac{\texte^{(\beta-1)\lambda^2|B_L|}}
{\beta^{\frac12(|\T_L|-1)}}\biggr]^{\frac1{4|\T_L|}}.
\end{equation}
Choosing~$\beta=\lambda^{-2}$, letting~$L\to\infty$ and~$p\uparrow\pt$, we thus conclude
\begin{equation}
\tilde\mu_\dis(\kappa_b\eta_b^2\le\lambda^2)\le c_1\lambda^{\ffrac14}.
\end{equation}
Noting that~$\tilde\mu_\dis(\kappa_b\eta_b^2\le\lambda^2)\ge\tilde\mu_\dis(\kappaD\eta_b^2\le\lambda^2)-\epsilon$, the bound \eqref{E:ener-bd2} is also proved.
\end{proofsect}

\section{Duality arguments}
\label{S:duality}\noindent
The goal of this section is to prove Theorem~\ref{T:dual}. For that we will establish an interesting duality that relates the model with parameter~$p$ to the same model with parameter~$1-p$.

\subsection{Preliminary considerations}
\noindent
The duality relation that our model \eqref{V-double} satisfies 
boils down, more or less, to an algebraic fact that the plaquette condition~\eqref{E:loop}, represented by the delta function~$\delta(\eta_{b_1}+\eta_{b_2}-\eta_{b_3}-\eta_{b_4})$, can formally be written as
\begin{equation}
\label{E:5.1}
\delta(\eta_{b_1}+\eta_{b_2}-\eta_{b_3}-\eta_{b_4})
=\int\frac{\textd\phi^\star}{2\pi}\texte^{\texti\phi^\star(\eta_{b_1}+\eta_{b_2}-\eta_{b_3}-\eta_{b_4})}.
\end{equation}
We interpret the variable~$\phi^\star$ 
as the \emph{dual field} that is associated with the plaquette $(\eta_{b_1},\eta_{b_2},\eta_{b_3},\eta_{b_4})$. As it turns out (see Theorem~\ref{T:dualita}), by integrating the~$\eta$'s with the~$\phi^\star$'s fixed a gradient measure is produced whose interaction is
the same as for the~$\eta$'s, except that the~$\kappa_b$'s get replaced by~$1/\kappa_b$'s. This means that if we assume that
\begin{equation}
\label{E:kO=kD}
\kappaO\kappaD=1,
\end{equation}
which is permissible in light of the remarks at the beginning of Section~\ref{sec2.2}, then the duality 
simply exchanges~$\kappaO$ and~$\kappaD$! We will assume that \eqref{E:kO=kD} holds throughout this entire section.

The aforementioned transformation works nicely for the plaquette conditions
which guarantee that the $\eta$'s can \emph{locally} be integrated back to the $\phi$'s. 
However, in two-dimensional torus geometry, two additional global constraints are also required to ensure the \emph{global} correspondence between the gradients~$\eta$ and the~$\phi$'s. These constraints, which are by definition built into the \emph{a priori} measure~$\nu_L$ from Sect.~\ref{S:model}, do not transform as nicely as the local plaquette conditions. To capture these subtleties, we will now define another \emph{a priori} measure that differs from~$\nu_L$ in that it disregards these global constraints.

Consider the linear subspace~$\XX^\star_L\supset \XX_L$ of~$\R^{\B_L}$ that is characterized by the equations $\eta_{b_1}+\eta_{b_2}-\eta_{b_3}-\eta_{b_4}=0$ for each plaquette $(b_1,b_2,b_3,b_4)$.
This space inherits the Euclidean metric from~$\R^{\B_L}$; we define~$\nu_L^\star$ as the corresponding Lebesgue measure on~$\XX_L^\star$ scaled by a constant~$C_L$ which will be determined momentarily. In order to make the link with~$\nu_L$, we define
\begin{equation}
\label{vert-hor}
\etavert=\sum_{x\in\T_L}\eta_{x+\hate_1}
\quad\text{and}\quad
\etahor=\sum_{x\in\T_L}\eta_{x+\hate_2}.
\end{equation}
Clearly,
\begin{equation}
\label{XL-XL*}
\XX_L=\bigl\{ \eta\in \XX^\star_L\colon\etavert=0,\etahor=0\bigr\}.
\end{equation}
Consider also the projection
$\Pi_L\colon\XX_L^\star\to\XX_L$ which is
defined, for any configuration $\eta\in\R^{\B_L}$, by
\begin{equation}
(\Pi_L\eta)_b=\begin{cases}
\eta_b-\frac1{L^2}\etavert,\qquad&\text{if }b\in\Bvert_L,
\\
\eta_b-\frac1{L^2}\etahor,\qquad&\text{if }b\in\Bhor_L.
\end{cases}
\end{equation}
Then we have:

\begin{mylemma}
There exist constants~$C_L$ such that, in the sense of distributions,
\begin{equation}
\label{E:nu-rep}
\nu^\star_L(\textd\eta)=\biggl(\,L^2\int_\R\,\textd\theta\prod_{(b_1,b_2,b_3,b_4)}
\!\!\!\delta(\eta_{b_1}+\eta_{b_2}-\eta_{b_3}-\eta_{b_4}-\theta)\,\biggr)\prod_{b\in\T_L}\textd\eta_b.
\end{equation} 
Moreover, we have
\begin{equation}
\label{nuL-rep}
\nu_L(\textd\eta)=\nu^\star_L(\textd\eta)\delta(\etahor)\delta(\etavert)
\end{equation}
and
\begin{equation}
\label{nu-nu-reverse}
\nu^\star_L(\textd\eta)=
\nu_L \circ \Pi_L(\textd \eta)\,\lambda_{\Pi_L}(\textd\eta).
\end{equation}
Here, $\lambda_{\Pi_L}(\textd\eta)$ is a multiple of the Lebesgue measure on the two-dimensional space 
$\Pi_L^{-1}(0)\cap \XX^\star_L$, which can be formally identified with~$\textd\etahor\textd\etavert$.
\end{mylemma}

\begin{proofsect}{Proof}
We begin with~\eqref{E:nu-rep}. Consider the orthogonal decomposition~$\R^{|\B_L|}=\XX_L^\star\oplus(\XX_L^\star)^\perp$. Clearly,~$\dim \XX_L^\star=L^2+1$. Choosing an orthonormal basis~$w_1,\dots,w_n$ in~$(\XX_L^\star)^\perp$ (where~$n=\dim(\XX_L^\star)^\perp=L^2-1$) the measure~$\nu_L^\star$ can be written as
\begin{equation}
\nu_L^\star(\textd\eta)=C_L\Bigl(\,\prod_{j=1}^n\delta(w_j\cdot\eta)\Bigr)\prod_{b\in\B_L}\textd\eta_b.
\end{equation}
Let~$\ell_\pi$ denote the vectors in~$\R^{|\B_L|}$ such that if~$\pi=(b_1,b_2,b_3,b_4)$ then~$\ell_\pi\cdot\eta=\eta_{b_1}+\eta_{b_2}-\eta_{b_3}-\eta_{b_4}$. Then~$\ell_\pi\in(\XX_L^\star)^\perp$ with all but one of these vectors linearly independent. This means  that we can replace the linear functionals~$\eta\mapsto w_j\cdot\eta$ by the plaquette conditions. Fixing a particular plaquette,~$\pi_0$, we find that
\begin{equation}
\label{E:5.10a}
\nu^\star_L(\textd\eta)=\biggl(\,\prod_{\pi\ne\pi_0}\delta(\ell_\pi\cdot\eta)\,\biggr)\prod_{b\in\T_L}\textd\eta_b
\end{equation}
provided that
\begin{equation}
C_L=\bigl|\det(w_j\cdot\ell_\pi)\bigr|=\sqrt{\det(\ell_{\pi}\cdot\ell_{\pi'})}.
\end{equation}
The expression \eqref{E:5.10a} is now easily checked to be equivalent to \eqref{E:nu-rep}: Applying the constraints from the plaquettes distinct from~$\pi$, we find that $\eta_{b_1}+\eta_{b_2}-\eta_{b_3}-\eta_{b_4}=(1-L^2)\theta$.
The corresponding $\delta$-function becomes $\delta(L^2\theta)$, and so we can set~$\theta=0$ in the remaining $\delta$-functions. Integration over~$\theta$ yields an overall multiplier $\int_\R\delta(L^2\theta)\textd\theta=1/L^2$.

In order to prove \eqref{nuL-rep}, pick a subtree~$\TT$ of~$\T_L$ as follows:~$\TT$ contains the horizontal bonds in $\{b_1+\ell\hate_1\colon\ell=0,\dots,L-2\}$ and the vertical bonds in $\{b_2+\ell\hate_1+m\hate_2\colon\ell,m=0,\dots,L-2\}$. As is easy to check,~$\TT$ is a spanning tree.
Denoting by~$\gamma_L$ the measure on the right-hand side of \eqref{nuL-rep} pick a bounded, continuous function~$f\colon\R^{|\B_L|}\to\R$ with bounded support and consider the integral~$\int f(\eta)\gamma_L(\textd\eta)$. The complement of~$\TT$ contains exactly~$L^2+1$ edges and there are as many $\delta$-functions in \eqref{E:5.10a} and \eqref{nuL-rep}, in which all~$\eta_b$, $b\not\in\TT$, appear with coefficient~$\pm1$. We may thus resolve these constraints and substitute for all~$\{\eta_b\colon b\not\in\TT\}$ into~$f$---call the result of this substitution~$\tilde f(\eta)$. Then we can integrate all of these variables which reduces our attention to the integral~$\int\tilde f(\eta)\prod_{b\in\TT}\textd\eta_b$. 

As is easy to check, the transformation~$\eta_b=\phi_y-\phi_x$ for~$b=(x,y)$ with the convention $\phi_0=0$ turns the measure~$\prod_{b\in\TT}\textd\eta_b$ into~$\delta(\textd\phi_0)\prod_{x\ne0}\textd\phi_x$ and makes~$\tilde f(\eta)$ into~$f(\nabla\phi)$. We have thus deduced
\begin{equation}
\int_{\R^{|\B_L|}} f(\eta)\gamma_L(\textd\eta)=
\int_{\R^{|\TT|}}\tilde f(\eta)\prod_{b\in\TT}\textd\eta_b=
\int_{\R^{|\T_L|}} f(\nabla\phi)\,\delta(\textd\phi_0)\prod_{x\ne0}\textd\phi_x.
\end{equation}
From here we get \eqref{nuL-rep} by noting that the latter integral can also be written as $\int f(\eta)\nu_L(\textd\eta)$.

To derive~$\nu_L^\star=\nu_L\circ\Pi_L\textd\etahor\textd\etavert$, let us write $\XX_L^\star=\XX_L\oplus(\Pi_L^{-1}(0)\cap\XX_L^\star)$. Since~$\nu_L^\star$ is the $C_L$-multiple of the Lebesgue measure on~$\XX_L^\star$ and since~$\etahor$ and~$\etavert$ represent orthogonal coordinates in $\Pi_L^{-1}(0)\cap\XX_L^\star$, we have
\begin{equation}
\nu_L^\star(\textd\eta)=C_L\lambda_{\XX_L}\circ\Pi_L(\textd\eta)\,L^{-2}\textd\etahor\,L^{-2}\textd\etavert,
\end{equation}
where~$\lambda_{\XX_L}$ is the Lebesgue measure on~$\XX_L$. Plugging into \eqref{nuL-rep} we find that~$\nu_L
=C_L L^{-4}\lambda_{\XX_L}$ which in turn implies \eqref{nu-nu-reverse}.
\end{proofsect}

\begin{myremark}
It is of some interest to note that the measure~$\nu_L^\star$ is also reflection positive for direct reflections. One proof of this fact goes by replacing the~$\delta$-functions in \eqref{E:nu-rep} by Gaussian kernels and noting that the linear term in~$\theta$ (in the exponent) exactly cancels. The status of reflection positivity for the diagonal reflections is unclear.
\end{myremark}

\subsection{Duality for inhomogeneous Gaussian measures}
Now we can state the principal duality relation. For that let~$\T_L^\star$ denote the dual torus which is simply a copy of~$\T_L$ shifted by half lattice spacing in each direction. Let~$\B_L^\star$ denote the set of dual edges. We will adopt the convention that if~$b$ is a direct edge, then its dual---i.e., the unique edge in~$\B_L^\star$ that cuts through~$b$---will be denoted by~$b^\star$. Then we have:

\begin{theorem}
\label{T:dualita}
Given two collections~$(\kappa_b)_{b\in\B_L}$ and~$(\kappa^\star_b)_{b\in\B_L}$ of positive weights on~$\B_L$, consider the partition functions
\begin{equation}
\label{E:Z-deff}
Z_{L,(\kappa_b)}=\int\nu_L(\textd\eta)\,\exp\Bigl\{-\frac12\sum_{b\in\B_L}\kappa_b\eta_b^2\Bigr\}
\end{equation}
and
\begin{equation}
\label{E:Zstar-deff}
Z_{L,(\kappa^\star_b)}^\star=\int\nu_L^\star(\textd\eta)\,\exp\Bigl\{-\frac12\sum_{b\in\B_L}\kappa^\star_b\eta_b^2\Bigr\}.
\end{equation}
If $(\kappa_b)_{b\in\B_L}$ and~$(\kappa_b^\star)_{b\in\B_L}$ are dual in the sense that
\begin{equation}
\label{kk-star}
\kappa^\star_{b^\star}=\frac1{\kappa_b},\qquad b\in\B_L,
\end{equation}
then
\begin{equation}
\label{E:Z-rep}
Z^\star_{L,(\kappa^\star_b)}=2\pi L^2\biggl[\,\prod_{b\in\B_L}\sqrt{\kappa_b}\biggr]Z_{L,(\kappa_b)}.
\end{equation}
\end{theorem}

\begin{proofsect}{Proof}
We will cast the partition function~$Z^\star_{L,(\kappa^\star_b)}$ into the form on the right-hand side of \eqref{E:Z-rep}. Let us regard this partition function as defined on the dual torus~$\T_L^\star$. The proof commences by rewriting the definition \eqref{E:nu-rep} with the help of \eqref{E:5.1} as
\begin{equation}
\label{nu-rewrite}
\nu^\star_L(\textd\eta)=\biggl(L^2\int_\R\textd\theta\int_{\R^{\T_L}}\!\prod_{x\in\tilde\T_L}\frac{\textd\phi_x}{2\pi}
\exp\Bigl\{\,\texti\sum_{x\in\T_L}\phi_x(\eta_{\pi^\star(x)}-\theta)\Bigr\}\biggr)
\prod_{b^\star\in\B_L^\star}\textd\eta_{b^\star},
\end{equation}
where $\eta_{\pi^\star(x)}$ is the plaquette curl for the dual plaquette~$\pi^\star(x)$ with the center at~$x$.
Rearranging terms and multiplying by the exponential (Gaussian) weight from \eqref{E:Zstar-deff}, we are thus supposed to integrate the function
\begin{equation}
L^2\int_\R\textd\theta\int_{\R^{\T_L}}\!\prod_{x\in\T_L}\frac{\textd\phi_x}{2\pi}
\exp\Bigl\{-\frac12\sum_{b^\star\in\B_L^\star}(\kappa^\star_{b^\star}\eta_{b^\star}^2-2\texti\eta_{b^\star}\nabla_b\phi)-\texti\theta\sum_{x\in\T_L}\phi_x\Bigr\}
\end{equation}
against the (unconstrained) Lebesgue measure~$\prod_{b\in\B_L}\textd\eta_b$. Here~$\nabla_b\phi=\phi_y-\phi_x$ if~$b=(x,y)$ is dual to the bond~$b^\star$. Completing the squares and integrating over the~$\eta$'s produces the function
\begin{equation}
L^2\biggl[\,\prod_{b^\star\in\B_L^\star}\frac 1{\sqrt{\kappa^\star_{b^\star}}}\biggr]\int_{\R}\textd\theta
\int_{\R^{\T_L}}\prod_{x\in\T_L}\textd\phi_x\,
\exp\biggl\{-\frac12\sum_{b\in\B_L}\frac1{\kappa^\star_{b^\star}}(\nabla_b\phi)^2-\texti\theta\sum_{x\in\T_L}\phi_x\biggr\}.
\end{equation}
Invoking~\eqref{kk-star}, we can replace all~$1/\kappa^\star_{b^\star}$ by~$\kappa_b$.
The integral over~$\theta$ then yields~$2\pi$ times the $\delta$-function of~$\sum_x\phi_x$ which---%
by the substitution $\phi_x\mapsto\phi_x+\frac1{|\T_L|} \phi_0$ that 
has no effect on
the rest of the integral---can be converted to~$2\pi\delta(\phi_0)$. Invoking 
the definition \eqref{E:nuL} of~$\nu_L$, this leads to the partition function \eqref{E:Z-deff}.
\end{proofsect}

\begin{myremark}
\label{remark:mira-dualita}
Let~$Q_{L,p}$ be the extended gradient Gibbs measure~$Q_L$ for $\varrho=p\delta_{\kappaO}+(1-p)\delta_{\kappaD}$
with parameter~$p$ and let~$Q_{L,p}^\star$ be the corresponding measure with the \emph{a priori} measure~$\nu_L$ replaced by~$\nu_L^\star$. Then the above duality shows that the law of~$(\kappa_b)$ governed by~$Q_{L,p}$ is the same as the law of its dual~$(\kappa_b^\star)$---defined via \eqref{kk-star}---in measure~$Q_{L,p_\star}^\star$, once~$p$ and $p_\star$ are related by
\begin{equation}
\label{pp-dual}
\frac p{1-p}\,\frac{p_\star}{1-p_\star}=\sqrt{\frac{\kappaD}{\kappaO}}.
\end{equation}
Indeed, the probability in measure~$Q_{L,p_\star}^\star$ of seeing the configuration~$(\kappa_b^\star)$ with~$\NO^\star$ ordered bonds and~$\ND^\star$ disordered bonds is proportional to~$p_\star^{\NO^\star}(1-p_\star)^{\ND^\star}Z_{L,(\kappa^\star_b)}^\star$.
Considering the dual configuration~$(\kappa_b)$ and letting~$\ND=\NO^\star$ denote the number of disordered bonds and~$\NO=\ND^\star$ the number of ordered bonds in~$(\kappa_b)$, we thus have
\begin{equation}
\label{dual-prop}
p_\star^{\NO^\star}(1-p_\star)^{\ND^\star}Z_{L,(\kappa^\star_b)}^\star
=2\pi L^2\bigl(p_\star\sqrt{\kappaO}\,\bigr)^{\ND}\bigl((1-p_\star)\sqrt{\kappaD}\,\bigr)^{\NO}Z_{L,(\kappa_b)}.
\end{equation}
For~$p$ and~$p_\star$ related as in \eqref{pp-dual}, the right-hand side is proportional to the probability of~$(\kappa_b)$ in measure~$Q_{L,p}$. 
We believe that the difference between the two measures disappears in the limit~$L\to\infty$ and so the $\kappa$-marginals of the states~$\tilde\mu_\ord$ and~$\tilde\mu_\dis$ at~$\pt$ can be considered 
to be 
dual to each other. However, we will not pursue this detail at any level of rigor.
\end{myremark}

\subsection{
Computing
the transition point}
\label{sec5.3}\noindent
In order to use effectively the duality relation from Theorem~\ref{T:dualita}, we have to show that the difference in the \emph{a priori} measure can be neglected. We will do this by showing that both partition functions lead to the same free energy. This is somewhat subtle due to the presence (and absence) of various constraints, so we will carry out the proof in detail.

\begin{myproposition}
\label{prop-FE}
Let~$\varrho_L(\textd\kappa)=\prod_{b\in\B_L}[p\delta_{\kappaO}(\textd\kappa_b)+(1-p)\delta_{\kappaD}(\textd\kappa_b)]$ and recall that~$Z_{L,V}$ denotes the integral of~$Z_{L,(\kappa_b)}$ with respect to~$\varrho_L$. Similarly, let~$Z_{L,V}^\star$ denote the integral of~$Z_{L,(\kappa_b)}^\star$ with respect to~$\varrho_L$. Then (the following limits exist 
as $L\to\infty$ and)
\begin{equation}
\label{limit-same}
\lim_{L\to\infty}\frac1{|\T_L|}\log Z_{L,V}=\lim_{L\to\infty}\frac1{|\T_L|}\log Z_{L,V}^\star
\end{equation}
for all~$p\in[0,1]$.
\end{myproposition}

Before we commence with the proof, let us establish the following variance bounds for homogeneous Gaussian measures relative to the \emph{a priori} measure~$\nu_L$ and~$\nu^\star_L$:

\begin{mylemma}
\label{L:var-bd}
Let~$\mu_L$ 
be the (standard) Gaussian gradient measure
\begin{equation}
\label{E:5.23a}
\mu_L(\textd\eta)\propto\exp\Bigl\{-\frac12\sum_{b\in\B_L}\eta_b^2\Bigr\}\nu_L(\textd\eta)
\end{equation}
and 
$\mu_L^\star$ be the measure obtained by replacing~$\nu_L$ by~$\nu^\star_L$. For~$m=1,\dots,L$, let
\begin{equation}
Y_m=\sum_{\ell=0}^{m-1}\eta_{(\ell\hate_1,(\ell+1)\hate_1)}.
\end{equation}
There exists an absolute constant~$c_3>0$ such that for all~$L\ge1$ and all~$m=1,\dots,L$,
\begin{equation}
\Var_{\mu_L}(Y_m)\le\Var_{\mu^\star_L}(Y_m)\le c_3(1+\log m).
\end{equation}
\end{mylemma}

\begin{proofsect}{Proof}
In measure~$\mu_L$, we can reintroduce back the fields~$(\phi_x)$ and~$Y_m$ then equals $\phi_{m\hate_1}$. Discrete Fourier transform implies that
\begin{equation}
\Var_{\mu_L}(\phi_x)=\frac1{L^2}\sum_{\bk\in\tilde\T_L}\frac{|1-\texte^{\texti x\cdot\bk}|^2}{\widehat D(\bk)},
\end{equation}
where $\tilde\T_L$ is the reciprocal torus and $\widehat D(\bk)=|1-\texte^{\texti k_1}|^2+|1-\texte^{\texti k_2}|^2$ is the discrete (torus) Laplacian. Simple estimates show that the sum is 
bounded by a constant times~$1+\log|x|$, uniformly in~$L$. 
Hence,~$\Var_{\mu_L}(Y_m)\le\tilde c_3(1+\log m)$ for some absolute constant~$\tilde c_3$.

As for the other measure, we recall the definitions \eqref{vert-hor} and use these to write~$\eta_b=\nabla_b\phi+
\frac1{L^2}
\etahor$ if~$b$ is horizontal (and $\eta_b=\nabla_b\phi+
\frac1{L^2}
\etavert$ if~$b$ is vertical). The fact that the Gaussian field is homogeneous implies---via \eqref{nu-nu-reverse}---that the fields~$(\phi_x)$ and the variables~$\etavert$ and~$\etahor$ are independent with~$(\phi_x)$ distributed according to~$\mu_L$ and~$\etavert$ and~$\etahor$ Gaussian with mean zero and 
variance~$\ffrac{L^2}2$.
In this case~$Y_m=\phi_{m\hate_1}+\frac m{L^2}\etahor$ and so we get
\begin{equation}
\Var_{\mu^\star_L}(Y_m)=\Var_{\mu_L}(Y_m)+
\frac{m^2}{2L^2}.
\end{equation}
But~$m\le L$ and so the correction is bounded for all~$L$.
\end{proofsect}

\begin{proofsect}{Proof of Proposition~\ref{prop-FE}}
The proof follows the expected line: To compensate for the lack of obvious subadditivity of the torus partition function, we will first relate the periodic boundary condition to a ``fixed'' boundary condition. Then we will establish subadditivity---and hence the existence of the free energy---for the latter boundary  condition.

Fix~$M>0$ and consider the partition function~$Z_{L,V}^{(M)}$ defined as follows. Let~$\Lambda_L$ be a box of $L\times L$ sites 
and consider the set~$\B^0_L$ of edges with \emph{both} ends in~$\Lambda_L$.
Let $\nu_L^{(M)}(\textd\eta)$ be as in \eqref{E:nuL} subject to the restriction that~$|\phi_x|\le M$ for all~$x$ on the \emph{internal} boundary of~$\Lambda_L$. 
Let
\begin{equation}
Z_{L,V}^{(M)}=\int \varrho_L(\textd\kappa) \int\nu_L^{(M)}(\textd\eta) \exp\Bigl\{-\frac12\sum_{b\in\B^0_L}\kappa_b\eta_b^2\Bigr\}.
\end{equation}
We will now provide upper and lower bounds between the partition functions~$Z_{L,V}$ (resp.~$Z_{L,V}^\star$) and~$Z_{L,V}^{(M)}$, for a well defined range of values of~$M$.

Comparing explicit expressions for~$Z_{L,V}$ and $Z_{L,V}^{(M)}$ and using~$\kappa_b\le\kappaO$, we get
\begin{equation}
\label{E:5.28}
Z_{L,V}\ge Z_{L,V}^{(M)}\exp\bigl\{-\tfrac12\kappaO(2M)^2(2L)\bigr\}.
\end{equation}
To derive an opposite inequality, note that for~$\kappa_b\ge\kappaD$ we get that $\Var_{Q_{L,(\kappa_b)}}(\phi_x)\ge\Var_{\mu_L}(\phi_x)/\kappaD$, where~$\mu_L$ is as in~\eqref{E:5.23a}. Invoking one more time the Gaussian identity $E(\texte^X)=\texte^{EX
+\frac12\Var(X)}$ in conjunction with Lemma~\ref{L:var-bd}, yields
\begin{equation}
Q_L(\phi_x\ge M)\le\exp\Bigl\{-\frac12\frac{M^2}{\kappaD c_3(1+\log L)}\Bigr\}.
\end{equation}
Hence, if~$M\gg\log L$ we have that with probability at least~$\ffrac12$ in measure~$Q_L$, \emph{all} variables~$\phi_x$ are in the interval~$[-M,M]$. Since the interaction that wraps~$\Lambda_L$ into the torus is of definite sign, it follows that
\begin{equation}
\label{E:5.27a}
Z_{L,V}\le 2 Z_{L,V}^{(M)}
\end{equation}
for all~$L$ and all~$M\gg\log L$.

Concerning the star-partition function, Lemma~\ref{L:var-bd} makes the proof of 
\eqref{E:5.27a}
exactly the same. As for the alternative of \eqref{E:5.28}, we invoke \eqref{nu-nu-reverse} and restrict all~$|\phi_x|$ on the internal boundary of~$\Lambda_L$ to values less than~$M$ and~$|\etahor|$ and~$|\etavert|$ to values less 
than~$ML^2$. Since~$|\eta_b|=|\nabla_b\phi+\frac1{L^2}\etavert|\le 2M+M=3M$ for every vertical bond that wraps~$\Lambda_L$ into the torus (and similarly for the horizontal bonds), we now get
\begin{equation}
\label{E:5.28'}
Z^\star_{L,V}\ge Z_{L,V}^{(M)}(2ML)^2\exp\bigl\{-\tfrac12\kappaO(3M)^2(2L)\bigr\},
\end{equation}
where the factor $(2ML)^2$ comes from the integration over~$\etavert$ and~$\etahor$. We conclude that, for~$\log L\ll M=o(\sqrt L)$, the partition functions~$Z_{L,V}$, $Z_{L,V}^\star$ and~$Z_{L,V}^{(M)}$ lead to the same free energy, provided at least one of these exists.

It remains to establish that the partition function~$Z_{L,V}^{(M)}$ is (approximately) submultiplicative for some choice of~$M=M_L$. 
Choose, e.g.,~$M_L=(\log L)^2$ and let~$p\ge1$ be an integer. If two neighbors have their~$\phi$'s between~$-M_L$ and~$M_L$, the energy  across the bond is at most~$\frac12\kappaO(4M_L)^2$. Splitting~$\Lambda_{pL}$ into~$p^2$ boxes of size~$L$, and restricting the~$\phi$'s to~$[-M_L,M_L]$ on the internal boundaries of these boxes, we thus get
\begin{equation}
Z_{pL,V}^{(M_{pL})}\ge[Z_{L,V}^{(M_L)}]^{p^2}\exp\bigl\{-\tfrac12\kappaO(2M_L)^2\,2(p-1)L\bigr\}.
\end{equation}
The exponent can be bounded below by $(pL)^{\ffrac32}-p^2L^{\ffrac32}=-(p^2-p^{\ffrac32})L^{\ffrac32}$ for~$L$ sufficiently large which implies that $p\mapsto[Z_{pL,V}^{(M_{pL})}\exp\{-(pL)^{\ffrac32}\}]^{1/(pL)^2}$ is increasing for all~$p\ge1$ and all~$L\gg1$. This proves the claim for limits along multiples of any fixed~$L$; to get the values ``in-between'' we just need to realize that, as before,~$Z_{L+k,V}^{(M_{L+k})}\ge Z_{L,V}^{(M_L)}\texte^{O(kLM_L^2)}$, for any fixed~$k$.
\end{proofsect}

Now we finally prove our claim concerning the value of the transitional~$p$:

\begin{proofsect}{Proof of Theorem~\ref{T:dual}}
Let~$Z_{L,V}^{(p)}$ denote the integral of~$Z_{L,(\kappa_b)}$ with respect to the \emph{a priori} measure 
$\varrho_L(\textd\kappa)$ with parameter $p$
and let~$Z_{L,V}^{\star,(p)}$ denote the analogous quantity for~$Z_{L,(\kappa_b)}^\star$. The arguments leading up to \eqref{dual-prop} then 
yield
\begin{equation}
Z_{L,V}^{\star,(p_\star)}=Z_{L,V}^{(p)}(2\pi L^2)\bigl(p_\star\sqrt{\kappaO}+(1-p_\star)\sqrt{\kappaD}\,\bigr)^{|\B_L|}
\end{equation}
whenever~$p_\star$ is dual to~$p$ in the sense of \eqref{pp-dual}. Thus, 
using~$F(p)$ to denote the limit in \eqref{limit-same} with the negative sign, we have
\begin{equation}
\label{FF-dual}
F(p_\star)=F(p)-2\log\bigl(p_\star\sqrt{\kappaO}+(1-p_\star)\sqrt{\kappaD}\,\bigr).
\end{equation}
Now, as a glance at the proof of Theorem~\ref{T:torus} reveals,
the value~$\pt$ is defined as the unique point where
the derivative of~$F(p)$, which at the continuity points of~$p\mapsto\chi(p)$ is simply $F'(p)=2\chi(p)-1$,
jumps from values near~$-1$ to values near~$+1$. 
Eq.~\eqref{FF-dual} then forces the jump to occur at the self-dual point~$p_\star=p$. 
In light of~\eqref{pp-dual}, this proves~\eqref{pt-eq}.
\end{proofsect}

\section*{Acknowledgments}
\noindent 
The research of M.B. was supported by the NSF grant~DMS-0505356 and that of R.K.~by the grants GA\v CR  201/03/0478,  MSM~0021620845, and the Max Planck Institute for Mathematics in the Sciences, Leipzig. The authors are grateful to Scott Sheffield for discussions that ultimately led to the consideration of the model~\eqref{V-double}, and for valuable advice how to establish the zero-tilt property of the coexisting~states. Discussions with Jean-Dominique Deuschel helped us understand the problems described in Remark~\ref{rem4.9}.

\end{document}